\input amstex
\documentstyle{amsppt}

\define\CC{\operatorname{\Bbb C}}
\define\VV{\Bbb V}
\define\UU{\Bbb U}
\define\NN{\Bbb N}

\define\jan{\goth s\goth l(n+1,\CC)}

\define\ban{\goth s\goth l(n-1,\CC)}
\define\bn{{{\goth s\goth o(2n+1,\CC)}}}
\define\jbn{{\goth s\goth o(2n+3,\CC)}}
\define\cn{{\goth s\goth p(2n,\CC)}}
\define\dn{{\goth s\goth o(2n,\CC)}}
\define\jdn{{\goth s\goth o(2n+2,\CC)}}

\define\semis{{(\gog_0^{\CC})_s}}
\define\ti{\tilde}

\define\gog{\goth g}
\define\goh{\goth h}
\define\RR{\operatorname{\Bbb R}}

\define\>{\rightarrow}
\define\<{\leftarrow}
\define\[{\lbrack}
\define\]{\rbrack}
\redefine\o{\circ}

\define\na{\nabla}
\define\al{\alpha}
\define\be{\beta}
\define\ga{\gamma}
\define\de{\delta}

\define\et{\eta}
\define\th{\theta}

\define\la{\lambda}
\define\rh{\rho}
\define\si{\sigma}
\define\ta{\tau}

\define\om{\omega}
\define\Ga{\Gamma}
\define\De{\Delta}

\define\La{\Lambda}

\define\Ph{\Phi}

\define\row#1#2#3{#1_{#2},\ldots,#1_{#3}}
\define\rowup#1#2#3{#1^{#2},\ldots,#1^{#3}}
\define\irow#1#2#3{{#1_{#2}\ldots#1_{#3}}}
\define\irowup#1#2#3{{#1^{#2}\ldots#1^{#3}}}

\define\x{\ctimes}
\define\sqr#1#2{{\vcenter{\vbox{\hrule height.#2pt
   \hbox{\vrule width.#2pt height#1pt \kern#1pt
         \vrule width.#2pt}
         \hrule height.#2pt}}}}
\define\sqre{\mathchoice\sqr{5.75}3\sqr{5.75}3\sqr{4.4}3\sqr{3}3}
\define\ctimes{\mathbin{\sqre\kern-6.5pt\times}}

\define\ox{\otimes}
\define\nmb#1#2{#2}      

\define\C{{\Cal C}}

\define\G{{\Cal G}}


\redefine\D{{\Cal D}}

\predefine\SS\S
\redefine\S{{\Cal S}}

\define\dd#1{\tfrac \partial{\partial #1}}
\define\field#1#2#3{{#1}^{#2}\frac\partial{\partial{#1}_{#3}}}

\define\pol#1 #2#3{L^{#1}_{\text{sym}}(\Bbb R^#2,\Bbb R^#3)}

\def\today{\ifcase\month\or
 January\or February\or March\or April\or May\or June\or
 July\or August\or September\or October\or November\or December\fi
 \space\number\day, \number\year}
\title Invariant operators on manifolds with almost Hermitian
symmetric structures, \\III. Standard operators.\endtitle
\author Andreas \v Cap, Jan Slov\'ak, Vladim\'ir Sou\v cek\endauthor
\abstract
This paper demonstrates the power of the calculus developed in the two
previous parts of the series for all real forms of the almost Hermitian
symmetric structures on smooth manifolds, including e.g\. conformal
Riemannian and almost quaternionic geometries. Exploiting some finite
dimensional representation theory of simple Lie algebras, we give explicit
formulae for distinguished invariant curved analogues of the standard
operators in terms of the linear connections belonging to the structures in
question, so in particular we prove their existence. Moreover, we prove that
these formulae for $k$th order standard operators, $k=1,2,\dots$, are
universal for all geometries in question.
\endabstract
\address
Institut f\"ur Mathematik, Universit\"at Wien, Strudlhofgasse 4,
1090 Wien, Austria\newline\indent
Department of Algebra and Geometry, Masaryk University in Brno,
Jan\'a\v ckovo n\' am\. 2a, 662~95 Brno, Czech
Republic\newline\indent
Mathematical Institute, Charles University, Sokolovsk\'a 83,
Praha, Czech Republic
\endaddress
\thanks The second and third authors supported by the GA\v CR,
grant Nr\. 201/96/0310
\endthanks
\rightheadtext{Invariant operators on manifolds with AHS structures III.}
\endtopmatter

\document \head\nmb0{1}. Introduction\endhead

As generally known, several geometries share surprisingly many properties
with the conformal Riemannian structures and projective structures. For
example the almost quaternionic ones. Following the old ideas by Cartan, and
some more recent development by Baston, Eastwood, Gindikin, Goncharov,
Ochiai, Tanaka, and others, we have started the project of building a good
calculus for all of them. This paper presents the first major application of
the technique developed so far for the so called AHS-structures in the first
two parts of this series, \cite{CSS1, CSS2}.

In \cite{F}, Fegan described all conformally invariant operators of the
first order on conformal Riemannian manifolds.  We use the invariant
differentiation with respect to Cartan connections developed in
\cite{CSS1}, together with some representation theory of simple Lie algebras,
in order to extend Fegan's methods to operators of all orders. This new
technique works for a wide class of geometries and, using the explicit
computations of the canonical Cartan connections in \cite{CSS2}, we obtain
formulae for all these invariant operators in terms of covariant derivatives
with respect to the linear connections belonging to the structures and their
curvatures. Moreover, a simple recursive procedure for the computation of
the correction terms for standard operators is described.

In such a way, the abstract indication of the existence of the standard
invariant linear differential operators on manifolds with almost Hermitian
symmetric structures given in \cite{B} is replaced by an explicit and
transparent construction, which provides even formulae in closed forms.
Surprisingly enough, these universal formulae do not depend on the
particular geometry at all.

In order to make the paper more self-contained, we have included a brief
review of some background from \cite{CSS1}. This concerns the short section
\nmb!{2} where we also fix the notation used in the sequel. The sections
\nmb!{3} through \nmb!{5} provide the necessary development in
representation theory. In order to address a wider audience among
differential geometers, we try to be quite detailed here. Section \nmb!{6}
gives the main existence result (Theorem \nmb!{6.5}) and the explicit
formulae are established in section \nmb!{7} (Theorems \nmb!{7.4} and
\nmb!{7.9}). Some technical points are postponed to two appendices.

\head\nmb0{2}. A calculus for Cartan connections\endhead

The aim of this section is to summarize for convenience
of the reader the main development from \cite{CSS1}. Full details
and proofs can be found there.

\subhead\nmb.{2.1} AHS structures\endsubhead
A basic datum distinguishing a particular AHS structure
is a real simple Lie group $G$ with the
Lie algebra $\gog$, which is $|1|$-graded, i.e\.
$$
\gog=\gog_{-1}\oplus\gog_0\oplus\gog_1
$$
with
$[\gog_i,\gog_j]\subset\gog_{i+j};\;\gog_j=\{0\},j\not=-1,0,1.$
There is a list of all simple real $|1|$-graded Lie algebras (see \cite{KN}).
Their complexification is a semisimple $|1|$-graded complex Lie
algebra. The classification of complex simple $|1|$-graded Lie
algebras corresponds
to the well known list of Hermitian symmetric spaces. The latter fact has
been the origin of the name A(lmost) H(ermitian) S(ymmetric) we use.

The subalgebras $\gog_{\pm 1}$ are commutative and dual to each other with
respect to the Killing form.
The algebra $\gog_0$ is  reductive with one-dimensional  center,
which is generated by the grading element $E,$ which is
characterized by the fact that each of the subalgebras $\gog_j$, $j=-1,0,1$,
ist the eigenspaces for the adjoint action of $E$ with eigenvalue $j$.
The semisimple  part $[\gog_0,\gog_0]$ of $\gog_0$ will
be denoted by $\gog_0^s$.

The subgroups $P,$ resp\. $P_1$
of $G$
correspond to the Lie algebra $\goth p=\gog_0\oplus\gog_1$,
resp\. $\gog_1$. The group $P_1$ is a normal subgroup of $P$
and the group  $G_0=P/P_1$ has the Lie algebra $\gog_0$. Let us mention
that we have used the letter $B$ instead of $P$ in \cite{CSS1}.

The typical and best understood example of AHS structures is
a conformal structure on a manifold $M$. A standard way
to define it is a reduction of the frame bundle of $M$
to the conformal group $G_0=CO(n,\RR)$.
A classical theorem going back to Cartan gives a
construction of a $P$-principal bundle $\G$ (where $P$
is a semidirect product of $G_0$ and $\RR^n$)
over $M$ and a uniquely defined Cartan connection $\om$
on $\G$. Such data were considered by Cartan as a curved analogue
of the flat model $G/P$ (an example of his `espaces
g\'en\'eralis\'es').
The characteristic properties of the Cartan connection $\om$
are a simple generalization of properties of the Maurer-Cartan
form $\om$ on $G/P$.

Following previous results by Tanaka, Ochiai, and Baston,
a simple and transparent principal bundle approach to a canonical
construction of the principal bundle $\G$ with structure group $P$
and of the Cartan connection $\om$
on $\G$ from the standard first order $G_0$-structure on $M$ was
described in \cite{CSS2}. We shall not need the construction here
and we shall start with $\G$ and  $\om$ as with a given
prescribed data, giving to $M$ the structure of an AHS manifold.

\subhead\nmb.{2.2} The Cartan connection and the
invariant differential\endsubhead
So we suppose that a $P$-principal bundle $\G$ on $M$ and
the Cartan connection $\om\in\Omega^1(\G,\gog)$ is given on $\G$
(for the definition and properties of the Cartan connections, see
\cite{CSS1}).

Any  Cartan connection defines an absolute
parallelism of $\G$ and for any vector space $\VV$, we can define the
{\it invariant differential}
$$
\na^{\om}:\C^{\infty}(\G,\VV)\rightarrow
\C^{\infty}(\G,\gog_{-1}^*\ox\VV)
$$
by $$
\na^{\om}s(u)(X)\equiv \na^{\om}_Xs(u):=
[\om^{-1}(X)s](u)
$$
where $\om^{-1}(X)$ is the constant vector field on $\G$ given by $X\in
\gog_{-1}$ and $\om$. Notice also $TM=\G\times_P \gog_{-1}$,
$T^*M=\G\times_P\gog_1$ in a canonical way.

If $\VV$ is a (finite dimensional) $P$-module,
than the space $\C^{\infty}(\G,\VV)^P$ of equivariant maps
is a 'frame form' of the space $\Ga(M,V)$ of smooth
sections of the associated vector bundle $V=\G\times_P\VV$.
We would like to use $\na^{\om}$ for a construction of
invariant differential operators. Unfortunately, the map
$\na ^{\om}s$, $s\in \C^{\infty}(\G,\VV)^P$, does not usually
belong to $\C^{\infty}(\G,\gog_{-1}^*\ox\VV)^P$,  it is not
the frame form of a section of a suitable associated vector
bundle over $M$. So $\na ^{\om}$ does not define directly
a differential operator on $M$.

A very useful procedure how to improve the situation is
to introduce a functorial way how to define a structure
of a $P$-module on the space
$$
J^1(\VV):=\VV\oplus (\gog_{-1}^*\ox \VV)
$$
in such a way that the map
$$
s\in\C^{\infty}(\G,\VV)^P\mapsto
(s,\na^{\om}s)\in\C^{\infty}(\G,J^1(\VV))^P
$$
has again  values in the space of equivariant maps.
The  $P$-module structure on $J^1(\VV)$ can be deduced
easily from the corresponding homogeneous case (where it is just
the representation inducing the homogeneous bundle  $J^1(V)$
of $1$-jets of sections of $V$).
Moreover, the Cartan connection $\om$ introduces the
natural identifications of the first jet prolongations of the
associated bundles $V=\G\times_P \VV$ with $\G\times_P J^1(\VV)$.

Consequently, any $P$-module homomorphism
$\Phi:J^1(\VV)\rightarrow \VV'$ induces a well
defined differential operator from the space of sections of the bundle $V$
to the space of sections of the bundle $V'$.
Due to the fact that the Cartan connection is uniquely defined
by the AHS structure, the  corresponding operator is invariant
with respect to any of the usual definitions of invariant
operators (details on relations between various possible definitions
of invariant operators can be found in \cite{Slo}).

The  situation most commonly considered is the case when $\VV$
and $\VV'$ are
irreducible $P$-modules. It means that $\VV$ (resp\. $\VV'$) are
 irreducible
$G_0$-modules with the trivial action of the nilpotent part of $P$.
In such a case, natural candidates for   $P$-homomorphisms $\Phi$
are projections from the space $\gog_{-1}^*\ox\VV$
(considered as an $\gog_0^s$-module) onto its irreducible factors,
extended by zero on the $\VV$ part of the module $J^1(\VV)$.
We shall see below that for any such projection, there is just
one specific value for the action of the grading element $E$ for
which the corresponding projection is a $P$-homomorphism
and that any invariant first order differential operator
on a manifold with a given AHS structure is obtained by this
construction.
For conformal structures, this was exactly the content of
the classification theorem obtained by Fegan in \cite{F}
(see \nmb!{7.2} below).

\subhead\nmb.{2.3} Iterated differentiation, semiholonomic jets\endsubhead
Iteratively, we can define the functor $\bar{J}^k(-)$
(the $k$-th semi-holonomic
prolongation) mapping any $P$-module $\VV$
to a submodule $\bar{J}^k(\VV)$ of the
$P$-module  $J^1(\bar{J}^{k-1}(\VV))$.
Considered as a $G_0$-module,  it looks like
$$
\bar{J}^k(\VV)=\VV\oplus(\gog^*_{-1}\ox\VV)\oplus...\oplus
(\ox^k(\gog^*_{-1})\ox\VV).
$$
As in the first order case,
the iterated invariant differential $(\na^{\om})^k$
defines the map
$$
j^k_{\om}:s\in\C^{\infty}(\G,\VV)^P\mapsto
(s,\na^{\om}s,\ldots,(\na^{\om})^k
s)\in\C^{\infty}(\G,\bar{J}^k(\VV))^P.
$$
Moreover, if $V=\G\times_P\VV$ is the bundle associated to $\VV$,
then its $k$th semi-holonomic jet prolongation
$\bar{J}^k(V)$ is the bundle associated to the
representation $\bar{J}^k(\VV)$. Thus construction of a large 
class of higher order invariant differential operators is now
possible as it was in the first order case: It is sufficient
to take any $P$-homomorphism from $\bar{J}^k(\VV)$ to
a $P$-module $\VV'$ and to compose it with the map $j^k_{\om}$.

The question to be answered is how to construct  such $P$-module
homomorphisms. If $\VV$ is an irreducible $P$-module, then it is easy
to find all $G_0$--module homomorphisms between the  corresponding
modules using representation theory. An explicit criterion showing
when such a $G_0$--homomorphisms is actually a $P$--module
homomorphism, was proved in \cite{CSS1} and will be used below to
prove existence results for invariant operators
(see \nmb!{5.2} for more details).

\subhead\nmb.{2.4}
Distinguished connections, the deformation tensor\endsubhead
Invariant operators are given as a composition of a
suitable $P$-homomorphism and the Cartan connection.
To express the result in standard terms (covariant derivatives,
curvature terms) and to find
explicit formulas for it, we need  more information.

Let us recall first the  relation between the original first order
structure $\G_0$ on $M$ (e.g\. a conformal one in the best known example)
and the $P$-principal bundle $\G$ constructed from it.
If $P_1$ is the Lie group corresponding to the Lie algebra
$\gog_1,$ then $\G_0\simeq \G/P_1$.
The value of the Cartan
connection $\om$ can be split with respect to the grading of
$\gog$ as $\om=\om_{-1}+\om_0+\om_1$.
For any $G_0$-equivariant
section $\si\:\G_0\to\G$ (which always exists),
the pullback $\si^*\om_0$ is a
principal connection on $\G_0.$ The space of all such connections
is an affine space modeled on the space of $1$-forms on $M$.
We have got in such a way a distinguished class of connections
on $M$ which are completely characterized by the requirements that
they have to belong to $\Cal G_0$, and their torsion has to coincide with the
$\gog_{-1}$-component of the curvature of $\om$.
In the conformal case, for example,  this class consists of all Weyl
geometries (thus contains all Levi-Civita
connections corresponding to any Riemannian metric chosen inside
the given conformal class, in particular).
The associated covariant derivatives
are standard tools used for description of differential
operators.



If $\om$ and $\tilde{\om}$ are
two Cartan connections which differ only in the
$\gog_1$-component, there exists an equivariant map
$\Ga\in{\Cal C}^{\infty}(\G,\gog_{-1}^*\ox\gog_1)$
such that $\tilde{\om}=\om-\Ga\circ\om_{-1}.$ The map $\Ga$
is the $P$-equivariant representation on $\G$ of a tensor
on $M,$ which is called the {\it deformation tensor}. In particular, once we
fix the Cartan connection $\om$ and the $G_0$-equivariant
section $\si\:\G_0\to \G$, there is the unique Cartan
connection $\tilde \om$ which is $\si$-related to the pullback
$\si^*(\om_{-1}+\om_0)$. This is the Cartan connection whose invariant
derivative $\na^{\tilde\om}$ is as close to the covariant derivative
$\na^{\si^*\om_0}$ as possible. The corresponding deformation tensor $\Ga$
then gives the full remaining comparison.
For conformal structures, this is just the well known  `rho--tensor'
having the following expression in terms of the Ricci curvature:
$$
\Ga_{ij}=\frac{-1}{m-2}\bigl(R_{ij} -
\frac{\de_{ij}}{2(m-1)}R\bigr),
$$
where $R_{ij}$ and $R$ are the $P$-equivariant
pull-backs of the Ricci tensor and
the scalar curvature to $\G$ and $m$ is the dimension of the manifold $M$.
Thus $\Ga$ is a generalization of the `rho--tensor'
to all AHS structures. Similar explicit formulae for these rho-tensors for
most AHS structures have been computed in \cite{CSS2}. 

Now, the value $\na^{\om}s$ of the invariant differential
on a section $s$ can be described in  more familiar terms,
using  $\na^{\ga}$ and the
deformation tensor $\Ga$ as follows.
The choice of $\si$ defines the trivialization
of the bundle $p\:\G\to\G_0$ expressed by the second coordinate
$\ta\:\G\to\gog_1$,
which can be characterized by the formula
$u=\si(p(u))\cdot\exp(\ta(u))$.
Let $\VV$ be an irreducible $P$-module,  $V=\G\times_P\VV\simeq
\G_0\times_{G_0}\VV$ the
corresponding associated vector bundle.
Sections $s\in\Ga(V)$ will be represented by means of  equivariant maps
$s\in{\Cal C}^{\infty}(\G_0,\VV)^{G_0}$
or equivalently as $p^*s\in{\Cal C}^{\infty}(\G,\VV)^P$.
Then we have for all $u\in P$, $X\in\gog_{-1}$
$$
\left(\na^{\om}(p^*s)(u)\right)(X)=
(p^*(\na^{\ga}s))(u)(X)
+[X,\ta(u)]\cdot ((p^*s)(u))
$$
where the bracket $[X,\ta(u)]\in \goth g_0$ acts on the element of the
$\gog_0$-module $\VV$.

All terms in the formula are $G_0$-equivariant, but only the
first one is also $P_1$-equivariant (i.e\.  constant along fibers
of $p$).
It is the map $\ta$ in the second term, which is not
$P_1$-equivariant (it varies when $u\in\G$ changes its position in the
fiber). This shows again that the invariant differential $\na^{\om}s$
is not $P$-equivariant even if $s$ itself is.
In many cases we can find a homomorphism $\Phi$ in such a way that
the term containing $\tau$ is killed by $\Phi$ and the
resulting composition is an invariant operator.

\subhead\nmb.{2.5} Correction terms and obstruction terms \endsubhead
To construct higher order invariant operators, we have to
use higher order iterations of the invariant differential.
To understand what is happening in higher orders, the second
order case is a representative example. It is possible
again to express $(\na^{\om})^2s$ using $\na^{\ga}$ and $\Ga$.
For any section $s\in{\Cal C}^{\infty}(\G_0,\VV)^{G_0},$ we have
$$\align
\left( (\nabla^\om)^2 (p^*s) \right)=\ &
 p^*( (\nabla^{\ga})^2s)+
D_0(\ga,\Ga)+D_1(\ga,\Ga,\ta)+D_2(\ga,\Ga,\ta)
\\\noalign{where}
D_0(\ga,\Ga)(u)(X,Y)=\ &[X, \Ga(u).Y]\cdot (p^*s(u));
\\D_1(\ga,\Ga,\ta)(u)(X,Y)=
\ &[X,\ta(u)]\cdot (p^*(\nabla^{\ga}_Ys))(u)+
\left([Y,\ta(u)]\cdot(p^*\nabla^{\ga}s)(u)\right)(X);
\\
D_2(\ga,\Ga,\ta)(u)(X,Y)=\ &
\left([Y,\ta(u)]\cdot([\_,\ta(u)]\cdot
(p^*s)(u))\right)(X) \\&-
\tfrac12[X,[\ta(u),[\ta(u),Y]]]\cdot (p^*s)(u),
\endalign
$$
and \ $\cdot$ \ denotes the appropriate action of an element
from $\gog_0$ on the space in question (either $\VV$ or
$\gog_{-1}^*\ox\VV)$.
The term $D_0$ is called the {\it correction term} and the terms
$D_i$, $i=1,2$, which are homogeneous of degree $i$ in $\ta$,
are called {\it obstruction terms}.

As for the first order case, the map
$(\nabla^\om)^2 (p^*s)$ is only
$G_0$-equivariant and, in general, not $P$-equivariant.
To define an invariant second order operator, it is necessary
to kill all obstruction terms by a suitable $G_0$-homomorphism.
If it is possible, then the leading term together with the
correction term gives an explicit formula for the corresponding
invariant operator (expressed already in standard language).

\subhead\nmb.{2.6} The algorithm for higher orders \endsubhead
In fact, it can be shown (see \cite{CSS1}) that vanishing of
$D_1(\ga,\Ga,\ta)$
implies vanishing of all higher order obstruction terms, so that
existence proofs can be simplified. The algebraic condition
discussed above is equivalent to vanishing of the sum of certain
terms linear in $\ta$, so that it is even more simple condition,
but it is only sufficient condition, not necessary one.

To have  an explicit algorithm for computation of the form of the
 correction terms, we need to take into account during the
inductive procedure all obstruction terms, not only the linear ones.
For that, we can use the  algorithm for recurrent
computation of the correction and obstruction terms, which was
proved in \cite{CSS1} (for more details see \nmb!{7.4}).
Using MAPLE, it was easy to implement this algorithm
and to compute explicitly the correction and obstruction terms
for low orders. The number of terms is growing enormously.
For the $6$th order, the full formula has 7184 terms and
the correction part itself has 328 terms.
We shall see later on that for standard operators studied below,
further essential simplification is possible and the final formula
will have only 10 summands.
To write down
on paper an explicit form of  invariant operators of higher orders
is too awkward. Nevertheless, we
shall see that for a broad class of operators, the algorithm
for the explicit form of the operator can be simplified
substantially and that the form of correction terms for standard
operators is remarkably stable
and universal, independently of the type of AHS structure and the 
representation $\VV$ considered (see section \nmb!{7}).

In the next sections, we shall use representation theory to show
how the theory explained above can be used for better
understanding of properties of standard invariant operators.

\head\nmb0{3}. $G_0$-homomorphisms\endhead

To construct invariant operators, we have to
learn how to construct $P$-homomor\-phisms
from $\bar{J}^k(\VV)$ to a $P$-module $\VV'$.
The first thing to do is to understand what are the possibilities
for $G_0$-homomorphisms. We shall concentrate on the
situation when $\VV$ is an irreducible $P$-module. This implies
that $\VV$ is an irreducible $G_0$-module and the nilpotent part
acts trivially. Representation theory offers enough tools
to classify  all $G_0$-homomorphisms in this case. Any such
homomorphism is equivalent to a projection of
$\bar{J}^k(\VV)$ onto one of its irreducible components and
a  decomposition of the tensor product
$\bar{J}^k(\VV)=(\ox^i\gog_{-1}^*)\ox\VV$
to irreducible components is a standard problem studied in
representation theory of semi-simple Lie groups. In this section,
we shall prove some additional facts needed for a construction
of $P$-homomorphisms and we shall deal with
a general complex semi-simple Lie
algebra $\gog$. Later on we shall use it for the semisimple part
$\gog_0^s=[\gog_0,\gog_0]$
of $\gog_0$.

\subhead\nmb.{3.1} Notation\endsubhead
Let us consider a complex semi-simple Lie algebra $\gog$ with a Cartan
subalgebra $\goh$, a set $\De^+$ of positive roots and its
subset $S=\{\al_1,\ldots,\al_n\}$ of simple roots.
Using the Killing form $(.,.)$, fundamental weights
$\pi_1,\ldots,\pi_n$ are defined by
$(\al_i^\vee,\pi_j)=\de_{ij}$,
where $\al_i^\vee =2\al_i/(\al_i,\al_i)$.

The (closed) dominant Weyl chamber
$\overline{\C}$ is given by linear combinations of fundamental
weights with nonnegative coefficients, let $\C$ denote
its interior.
Finite dimensional complex irreducible representations of
$\gog$ are characterized by their
highest weights $\la$, which lie in the weight lattice 
$\La^+=\{\sum\la_i\pi_i;\ \la_i\geq 0,\ \la_i\in\Bbb Z\}$.
The corresponding representation
will be denoted by $({\la},\VV_{\la})$ but the action
${\la}(X)v$, $X\in\gog$, $v\in \VV_{\la}$ will be often written simply
as $X\cdot v$, if the representation is clear from the context.
The set of all weights of $\VV$ will be denoted by
$\Pi(\VV).$

Any weight $\la\in\goh^*$ can be characterized by its
coefficients $\la_j=(\la,\al_j^\vee )$. In particular, the
simple  roots $\al_i$ have coefficients $a_{ij}=(\al_i,\al_j^\vee )$,
where $a_{ij}$ is the Cartan matrix of the Lie algebra $\gog$,
which is encoded into its Dynkin diagram.
Consequently, the reflection
$\si_i(\la)=\la-(\la,\al_i^\vee )\al_i$  with respect to a simple
root $\al_i$
changes coefficients $\la_j$ of $\la$ into
coefficients $\la_j-\la_i a_{ij}$.
Due to properties of the Cartan matrix,
the coefficient $\la_i$ changes to $-\la_i$
and (if no multiple edges of the Dynkin diagram are involved),
the coefficient $\la_i$ adds to neighboring coefficients $\la_j$
(for which $a_{ij}=-1)$.

The reflections $\si_i$ generate the Weyl group $W$. For 
$\rho=\sum_i\pi_i$, we shall denote  by $\cdot$ the affine
action of $W$ on weights defined by $w\cdot \la=
w(\la+\rho)-\rho$.

In our applications of the theory, we shall mostly need the case
of a simple Lie algebra $\gog$. The only exception will be the
Grassmannian case, where our Lie algebra $\gog$ will have two
simple parts $\gog^1\oplus\gog^2$. Note that in this case,
the Cartan subalgebra $\goh$ splits also into
$\goh^1\oplus\goh^2$, all weights can be written as couples
$\la=(\la^1,\la^2)$ and the representation $\VV_{\la}$ is the
tensor product $\VV_{\la^1}\ox\VV_{\la^2}$. The Killing form
splits as well: $(\la,\mu)=(\la^1,\mu^1)+(\la^2,\mu^2)$.
The Weyl group $W$ is the direct product $W_1\times W_2$ of the
Weyl groups of $\gog_1$ and $\gog_2$.

\subhead\nmb.{3.2} Klimyk's algorithm \endsubhead
There is a useful and explicit algorithm for the decomposition
of the tensor product of two irreducible representations
of a simple Lie algebra $\gog$
into irreducible components, based on the Klimyk formula
(see \cite{H}, Sec.24, Ex.9).

For any weight $\xi\in\goh^*,$ let $\{\xi\}$ denote the
dominant weight lying on the orbit of $\xi$ under the Weyl group.
If $\{\xi\}\in{\Cal C},$ then there is the unique $w\in W$
such that $\{\xi\}=w\xi$.
Let $t(\xi)$ be equal to the sign of $w$ in this case
and zero otherwise.

Suppose moreover that we know the list $\Pi(\mu)$ of all weights
of the irreducible representation $V_{\mu}$ with the highest
weight $\mu$, including their multiplicities $m_{\mu}(\nu)$, for 
$\nu\in\Pi(\mu)$.
Let $\VV_{\la}$ denote the irreducible representation of $\gog$
with the highest weight $\la$.
Then the Klimyk formula implies that
it is sufficient to go through the list $\Pi(\mu)$,
write a formal sum
$$\sum_{\nu\in\Pi(\mu)}m_{\mu}(\nu)t(\la+\rho+\nu)
\VV_{\{\la+\rho+\nu\}-\rho}
$$
of irreducible representations and to add together coefficients
at representations with the same highest weight.
The resulting coefficients are always non-negative
and give the multiplicity of the corresponding representation
in the decomposition.
Note that some cancelations happen often.

\subhead\nmb.{3.3} The decomposition of a tensor product of
representations
\endsubhead
There are certain facts known for a general case of a tensor
product of two irreducible representations
$\VV_{\la}$ and $\VV_{\mu}$ with highest weights $\la$ and $\mu$.
For example, the highest weight
$\xi$ of an irreducible piece  in the decomposition of the product
$\VV_{\la}\ox \VV_{\mu}$ has always form $\xi=\la+\nu,\nu\in\Pi(\mu)$
(see \cite{FH}, p.425). But in general,
we know nothing about its multiplicity,
it can be zero, one or bigger.

In the product $\VV_{\la}\ox \VV_{\mu},$ there is
always an irreducible piece with the highest weight $\la+\mu$
and it appears  with multiplicity one. This special
irreducible component is standardly denoted by
$\VV_{\la}\x \VV_{\mu}$, and called the {\it Cartan product}
of $\VV_{\la}$ and $\VV_{\mu}$.
If $e_{\la}$, resp\. $e_{\mu}$, are weight vectors for highest
weights $\la$, resp\. $\mu$, then $e_{\la}\ox e_{\mu}$
is a weight vector with the weight $\la+\mu.$
Consequently, $\x^k\VV\subset\odot^k\VV.$

The following general fact is much more difficult to verify.
The Parthasarathy--Rao--Varadarajan (PRV) conjecture proved
recently (see \cite{Ku}) claims that for any $w\in W$,
the module $\VV_{\{\la+w\mu\}}$ with the extremal weight
$\la+w\mu$ occurs in $\VV_{\la}\ox \VV_{\mu}$ with multiplicity
at least one.

In the case that one representation in a tensor product is in a
suitable sense small, we can say more about the decomposition. In
particular, there will be no multiplicities 
in the product for such cases. This is a substantial
information needed in applications below.
The simplest case is the following theorem.

\proclaim{Theorem}
Let $\mu$ be such that all weights $\nu\in\Pi(\mu)$ have multiplicity
one. Let us suppose moreover 
that the coefficients of all weights $\nu\in\Pi(\mu)$ with respect to
fundamental weights are $\ge-1$. Then for any $\la\in\La^+$, we have
$$
\VV_{\mu}\ox \VV_{\la}=\sum_{\ta\in A}\VV_{\ta}
$$
where  $A$ is the set of all weights
of the form $\ta=\la+\nu,\nu\in\Pi(\mu)$,
which belong to the dominant Weyl chamber $\overline{\C}$.
There are no multiplicities in the decomposition.
\endproclaim

\demo{Proof}
The coefficients in the decomposition of any weight $\la\in\La^+$
into fundamental weights are, by definition,  all nonnegative.
The weight $\rho$ has all coefficients equal to $1$.
Our assumptions above imply that for all weights $\nu\in\Pi(\mu)$,
the sum $\rho+\nu$ belongs to $\overline{\C}$,
hence $\la+\rho+\nu\in\overline{\C}$ as well.
So no action of elements $w\in W$ is needed,
$\{\la+\rho+\nu\}-\rho=\la+\nu$ for all $\nu\in\Pi(\VV_{\mu})$
and no cancelations or multiplicities in the
decomposition of the tensor product can occur.
The weight $\la+\nu$ appears in the decomposition (with nonzero
coefficient) if and only if 
$\la+\rho+\nu$ belongs to the interior $\C$ i.e\.
if and only if $\la+\nu\in\overline{\C}$.
\qed\enddemo


The theorem just proved will be sufficient in most cases needed
below. In two of them, we shall however need a case when
some of components of weights will be equal to $-2$.
We are going to prove the multiplicity one result
for this case under a suitable additional assumption. In some
particular cases (e.g\. in two cases needed below, see Appendix A),
it is possible
to describe the set $A$ in the decomposition more precisely,
but we shall not need to formulate such results in general.

\proclaim{Theorem'}
Suppose that $\mu$ is such that all weights
$\nu\in\Pi(\mu)$ have multiplicity one. Let us suppose moreover
that for all weights $\nu\in\Pi(\mu)$, $\nu=\sum_i\nu_i\pi_i,$
the following conditions are satisfied:
\roster
\item $\nu_i\geq -2$ for all $i$;
\item
there exists at most one index $i$ such that $\nu_i=-2$
and if it happens,
we suppose moreover that
for all $j\not= i$, $ \nu_j\geq 0$ and
$a_{ij}\geq -1$
(the last condition means that
the $i$th node of the corresponding Dynkin diagram
is not at the foot point of a double arrow).
\endroster
\noindent Then for any $\la\in\La^+,$ we have
$$
\VV_{\mu}\ox \VV_{\la}=\sum_{\ta\in A}V_{\ta}
$$
where  $A\subset(\{\la+\nu|\nu\in\Pi(\mu)\})\cap\overline{\C}$ is some
subset and there are no multiplicities in the decomposition.
\endproclaim

\demo{Proof}
For all weights $\nu$ with the  property
$\nu_j\geq -1$ for all $j$
we get as above that $\la+\nu+\rh\in\overline{\C}$,
 hence no reflections are needed and $V_{\la+\nu}$
appears in the formal sum coming from the Klimyk formula
if and only if $\la+\nu\in\overline{\C}$.

Let us consider a weight $\nu$  with the  property
that $\nu_i=-2$.
The assumptions of the theorem imply that
$(\la+\nu+\rho)_j\geq 1$, $j\not= i$, and
 $(\la+\nu+\rho)_i=\la_i-1$.
If $\la_i>0,$ then again $\la+\nu+\rho\in\overline{\C}$
and no reflection is needed.

If, however, $\la_i=0$ then
the weight $\la+\nu+\rho$ is not in $\overline{\C}$.
Let $w\in W$ is the simple reflection with respect to $i$th
simple root, then $(\la+\nu+\rho)_i=-1$ and
$(w\,(\la+\nu+\rho))_i=1$. For $j\not=i$ such that  $a_{ij}=0,$
 the coefficient $(\la+\nu+\rho)_j$ is not changed under the
reflection, hence is
nonnegative. If  $j\not=i$ such that  $a_{ij}=-1$, then
$(w\,(\la+\nu+\rho))_j=(\la+\nu+\rho)_j-1\geq \rho_j-1=0$,
hence also these coefficients are nonnegative.
Consequently, $w\,(\la+\nu+\rho)\in\overline{\C}$
and the irreducible representation $\VV_{w\,(\la+\nu+\rho) -\rho}$
will appear in Klimyk's formal sum with coefficient $-1$.

All terms in the formal sum coming from the weights $\nu$ with
the property
$\la+\nu+\rho\in\overline{C}$
are distinct and with multiplicity one.
All others are coming with the coefficients $-1$, hence they are necessarily
canceled by some of previous ones.
Hence all terms in the result have multiplicity one
and their highest weights are
contained in $\{\mu=\la+\nu,\nu\in\Pi(\mu)\}\cap\overline{\C}$.
\qed\enddemo

\subhead\nmb.{3.4} Multiple decompositions\endsubhead
We shall also have to understand
irreducible components of a more complicated
tensor product $(\ox^k\VV_{\mu})\ox\VV_{\la}$.
For $k>1$, there is no hope to get a multiplicity
one result as before.
As a consequence, only isotypic components of the product
will be unambiguously defined and the complete splitting
into irreducible components will depend on arbitrary choices.
We shall show now that the results of the previous paragraph
can be used for a
classification of the pieces in the decomposition
and for a construction of a distinguished decomposition
useful for more detailed computations in following sections.

Let $\gog$ is a semi-simple Lie algebra and $\VV_{\mu}$ its irreducible
representation having the following property:
For all $\la\in\La^+$, there exists
a set $A_{\la}$ such that $\VV_{\mu}\otimes\VV_{\la}=
\sum_{\la_1\in A_{\la}}\VV_{\la_1}$ and there are no multiplicities
in the decomposition.

Then the decomposition can be iterated as follows.
The product $\otimes^2(\VV_{\mu})\otimes \VV_{\la}=
\VV_{\mu}\otimes (\sum_{\la_1\in A_{\la}}\VV_{\la_1})$
can be again decomposed in the same way as
$$
\sum_{\la_1\in A_{\la}}\sum_{\la_2\in A_{\la_1}}\VV_{\la_2,\la_1},
$$
where the double index of $\VV_{\la_2,\la_1}$ indicates how this
particular component was obtained in the decomposition.
By repeating this process, it is clear that
the product $\otimes^k(\VV_{\mu})\otimes \VV_{\la}$
can be completely decomposed into irreducible components, each
one  being labeled by a sequence
$\underline{\la}=(\la_k,\la_{k-1},\ldots,\la_1,\la)$
 which records the way how
this  component was obtained through the process of
successive decompositions. The final highest weight $\la_k$ may
appear many times and its precise position in the
isotypic component is fixed by the whole sequence recording its
history.  Hence for a fixed $\la,$
 we shall define the set
$A_k(\la)$ of all such sequences, i.e.
$$
A_k(\la)=\{\underline{\la}=
(\la_k,\la_{k-1},\ldots,\la_1,\la_0)
\,|\,\la_0=\la,\,\la_j\in A_{\la_{j-1}};\;j=1,\ldots,k\}.
$$
Then
$$
\otimes^k(\VV_{\mu})\otimes \VV_{\la}=
\sum_{\underline{\la}\in A_k(\la)}\VV_{\underline{\la}}.
$$
Together with the final irreducible component
$\VV_{\underline{\la}}$, we shall use also for computations
all intermediate components
given by $\VV_{\underline{\la}_j},\underline{\la}_j=
(\la_j,\ldots,\la_0)$ in $\otimes^j(\VV_{\mu})\otimes \VV_{\la}$,
together with the corresponding invariant projections
$\pi_{\la_j}$.

There is one important question connected with such a
decomposition, namely to find a position of the above mentioned
components with respect to the splitting
of $\otimes^j(\VV_{\mu})\otimes \VV_{\la}$ into a direct sum of
$\odot^j(\VV_{\mu})\otimes \VV_{\la}$ and its invariant complement.
Such a knowledge would help to decide whether invariant operators
obtained by the projection to the corresponding components in the
decomposition will have nontrivial symbol or not.
We shall answer this question in the case  we need in the next
paragraph.

\subhead\nmb.{3.5} Multiplicity one components\endsubhead
There are special pieces in the
decomposition  of $\otimes^j(\VV_{\mu})\otimes \VV_{\la}$ which
 always appear with multiplicity one. Even more, we shall
be able to show  that they must be included in
$\x^j(\VV_{\mu})\otimes \VV_{\la},$ where
$\x$ denotes the Cartan product of irreducible
representations (cf\. \nmb!{3.3}),
hence their symbol will be nontrivial.

\proclaim{Theorem}
Let $\la,\mu\in\La^+$.
Let $\nu$ be an extremal weight of $\VV_{\mu}$
(i.e\. it belongs to the Weyl
orbit of the highest weight $\mu$).
 Let $k$ be a positive integer such that
$\la+k\nu$ is  dominant.

Then there is a unique  irreducible component in
$\ox^k(\VV_{\mu})\ox \VV_{\la}$ with
highest weight $\ta=\la+k\nu$.
Moreover, the component $\VV_{\ta}$ is contained in
$\x^k(\VV_{\mu})\ox \VV_{\la}$.
\endproclaim

\demo{Proof}
The product
$\ox^k(\VV_{\mu})\ox \VV_{\la}$
can be decomposed into the sum of $V_{\underline{\la}}$ as described
above. All these chains ${\underline{\la}}$
can be considered as piecewise linear paths in the
dominant Weyl chamber composed from the straight segment  with
directions
given by weights of $\VV_{\mu}.$ If we are going straight on
$k$ times in the same  direction given by an extremal weight of
$\VV_{\mu},$ no other path
can reach the same point $\ta=\la+k\nu$ (extremal weights have extremal
lengths). This implies the unicity of the
component.

To prove the existence, note that
the weight $k\nu$ is an extremal weight of $\x^k(\VV_{\mu}).$
Hence we can use the PRV conjecture
to show that  $\VV_{\ta} $ appears in the decomposition of
$\x^k(\VV_{\mu})\ox \VV_{\la}.$
\hfill\qed
\enddemo

\subhead\nmb.{3.6} Partial projections\endsubhead
Let us recall that we always have $\x^k(\VV)\subset
\odot^k(\VV)$
and that $\x^k(\VV)$ coincides with
$$
[\x^2(\VV)]\x[\x^{k-2}(\VV)]
\subset
[\x^2(\VV)]\ox[\x^{k-2}(\VV)].
$$
As a corollary we get

\proclaim{Lemma}
Denote by $\pi$ the projection of
$\ox^k(\VV)$
onto $\x^k(\VV).$
Suppose that $A$ is the invariant complement
of $\x^2(\VV)$ in $\ox^2(\VV)$
and $\pi_A$ is the corresponding projection.
Then $\left[A\ox(\ox^{k-2}(\VV)\right]\cap
\left[\x^k(\VV)\right]=\emptyset$, or equivalently
$$
\pi\o(\pi_A\ox Id^{\;k-2})=0.
$$
\endproclaim

\subhead\nmb.{3.7}\endsubhead
The results above will be applied below in the following special
case. Let $\gog=\gog_{-1}\oplus\gog\oplus\gog_1$
be a complex $|1|$-graded Lie algebra, cf\. \nmb!{2.1}. The space $\gog_1$
is an irreducible $\gog_0^s$-module which is `small' enough, i.e\.
it satisfies assumptions of one of the Theorems in \nmb!{3.3}.
To check it, it is necessary to inspect
algebras $\gog$ case by case. The list of them together
with details needed for the verification are collected in
Appendix A.

Consequently, for any irreducible $\gog_0$-module $\VV,$
the tensor product $\gog_1\ox\VV$ decomposes into irreducible
components without multiplicities and  results of \nmb!{3.5} and \nmb!{3.6}
can be used for decompositions of the product
$\ox^k(\gog_1)\ox\VV.$


\head\nmb0{4}. Casimir computations\endhead
\subhead\nmb.{4.1} Notation\endsubhead
For this section, we shall suppose that
$
\gog=\gog_{-1}\oplus\gog_0\oplus\gog_1
$
is  a complex $|1|$-graded simple Lie algebra.
In general, a choice of $|k|$-graded structure on a complex simple Lie
algebra $\gog$ is the same as a choice of its
parabolic subalgebra.
Any parabolic subalgebra is conjugated to a standard one
(i.e\. one containing a chosen Borel subalgebra $\goth b\subset
\gog$).
There is one to one correspondence between
standard parabolic subalgebras of $\gog$ and subsets of the
set $S$ of simple roots of $\gog$.

The $|1|$-graded structures on $\gog$
exist only for four classical series and for $E_6$ and $E_7$
cases and they are given by  certain one-point subsets of $S$ (Dynkin
diagrams with the corresponding simple root crossed are often used
to denote the chosen parabolic subalgebra).
We shall choose numbering of the set $S$ of simple roots
so that the first simple root $\al_0$ is the crossed one
(for more information on $|k|$-graded Lie algebras see \cite{BasE, Y}).

There is a unique grading element $E\in\gog_0$
satisfying $[E,X]=\ell X $ for $X\in\gog_{\ell},\,\ell=-1,0,1.$
A Cartan subalgebra $\goh\subset\gog$
can be chosen in such a way that
$E\in\goh,$ then $\goh\subset\gog_0.$
The set $\De^+$ of positive roots for $\gog$ can
be chosen so that all root spaces for positive roots
are included in $ \gog_0\oplus\gog_1$.

It is often useful to normalize an invariant form $(.,.)$ on $\gog$
by the requirement $(E,E)=1$ (see e.g\. \cite{BOO}).
For the Killing form, we have $B(E,E)=2\dim \gog_1$, hence
$(X,Y)=({2\dim\gog_1})^{-1}B(X,Y)$.
This normalized form $(.,.)$ induces
nondegenerate invariant bilinear forms on $\gog_{0}$ and
$
\gog_{-1}\oplus
\gog_{1}
$,
and it identifies $\gog_1$ and $\gog_{-1}$ as dual spaces.
Orthonormal bases and Casimir operators for $\gog_0$ will be computed using
this normalized form.

The algebra $\gog_0$ splits into $1$-dimensional center
$\goth a$
and a semisimple part $\gog_0^s=[\gog_0,\gog_0]$
which has $\goh_s=\goh\cap{\gog_0^s} $
as a Cartan subalgebra.
Then $\goh={\goth a}\oplus{\goh}_s$.
Irreducible  representations of
$\goth p=\gog_0\oplus\gog_1$ are
trivial on $\gog_1$. Every such representation is a
tensor product of a one-dimensional representation of $\goth a$
and an irreducible representation of $\gog^s_0$, which can be
characterized by its highest weight $\la\in(\goh^s)^*$.
For convenience, we shall consider $(\goh^s)^*$ as
a subset of $\goh^*$ of all elements, which restrict to zero on $\goth a$.
Representations  of $\goth a$ can be
characterized by a (generalized) conformal weight $w\in\CC.$
We shall say that a representation  $\UU$ of $\gog_0$ has a (generalized)
{\it conformal weight} $w$, if $E\cdot v=wv,\,v\in\UU$.
The cotangent spaces of our manifolds are associated to the
adjoint representation of $\gog_0$ on $\gog_1$, hence $1$-forms will
have (generalized) conformal weight $1$.
An irreducible representation of $\gog_0$ with a conformal
weight $w$ and highest weight $\la\in(\goh^s)^*$ will be
denoted by $\VV_{\la}(w)$.

Let $\{Y_a\}$, $a=0,1,\ldots$, be an orthonormal basis of $\gog_0$
with respect to the form $(.,.)$. We may choose it in
such a way that $Y_0=E\in\goth a$ and
$\{Y_{a'}\}$, $a'>0$ is an orthonormal basis for $\gog^s_0$.
For any representation $\VV$ of $\gog_0^s$, the Casimir operator
$C(\VV)$ is defined by
$C(\VV)=\sum_{a'>0}Y_{a'}\circ {Y}_{a'}$.
It is well known (see \cite{H}) that  if $\VV$ is an irreducible
representation with a highest weight $\la$, then
$$
C(\VV)=(\la,\la+2\rho); \quad\rho=1/2\sum_{\al\in\Delta^+(\gog_0^s)}
\al.
$$

As we have noticed already, our algebras $\gog_0^s$ are irreducible
in all cases except the $\goth s\goth l(n,\Bbb C)$ series, but even then
the formula $C(\VV_{\la})=(\la,\la+2\rho),\,\rho=(\rho_1,\rho_2)$
is still valid, see \nmb!{3.1} for the reasons.

\subhead\nmb.{4.2} Casimir computations\endsubhead
Suppose now that $X\in\gog_{-1}$, $Z\in\gog_1$
and let us consider an irreducible $\gog_0$-module
$\VV_{\la}(w),$ where $\la\in\goh^*$ is an integral dominant
weight for $\gog_0^s$ and $w\in\CC$.
In the description of iterated invariant differentials,
terms of type $[Z,X]\cdot s$, $s\in\VV_{\la}(w)$,
have appeared  very often (the $\cdot$ means here the action
of an element of $\gog_0$ under the representation characterized
by $\la$ and $w$),
(see \nmb!{2.4}). It is hence important to
understand them better.

Recall that we identify $\gog_1$ and $(\gog_{-1})^*$ using the scalar
product $(.,.)$. The term $[Z,X]\cdot s$ defines a map from
$\gog_1\ox\gog_{-1}\ox\VV_{\la}(w)$ into $\VV_{\la}(w),$
which can be interpreted also as a map
$\Phi\:\gog_1\ox\VV_{\la}(w)\to\gog_1\ox\VV_{\la}(w),$
defined by
$$
\Phi(Z\ox v)(X):={\la}([Z,X])v;\quad Z\in\gog_1,s\in\VV_{\la}(w),
X\in\gog_{-1}.
$$

Let  us choose bases $\{\et_{\al}\}$, resp\.
$\{\xi_{\al}\}$  of
$\gog_{-1}$, resp\. $\gog_1$, which are dual with respect
to the scalar product $(.,.)$. Due to
$$
[Z,X]\cdot s=\sum_{\al}[Z,(\et_{\al},X)\xi_{\al}]\cdot s=
\biggl(\sum_\al\et_{\al}\ox [Z,\xi_{\al}]\cdot s\biggr)(X),
$$
we get
$$
\Phi(Z\ox s)=\sum_\al\et_{\al}\ox [Z,\xi_{\al}]\cdot s.
$$

The map $\Phi$ is a $\gog_0$-homomorphism (by direct computation
or by the lemma below).
Let $\gog_1\ox \VV_{\la}(w)=\sum_{\mu}\VV_{\mu}(w+1)$
be a decomposition of the product of $\gog_0$-modules into
irreducible components and let
$\pi_{\la\mu}:\gog_1\ox\VV_{\la}(w)\>\VV_{\mu}(w+1)$
be the corresponding projections.  The $\gog_0$-homomorphism $\Phi$
acts as a multiple of the identity on each irreducible component,
i.e. there are constants ${\tilde c}_{\la\mu}\in{\Bbb R}$ such that
$\Phi=\sum_{\mu}\tilde{c}_{\la\mu}\pi_{\la\mu}$
and we are going to describe a formula expressing these constants
in terms of the weights $\la$ and $\mu$.

\proclaim{\nmb.{4.3} Lemma}
Let $\VV_{\la}(w)$  be an irreducible representation of
$\gog_0$ and
let
$\gog_1\ox \VV_{\la}=\sum_{\mu}\VV_{\mu}$
be a decomposition of the product into irreducible
$\gog_0^s$-modules.
Let
$\al$ be the highest weight of $\gog_1$ and let
$\rho$ be the half sum of positive roots for $\gog_0^s$. Then
for all $s\in\VV_{\la}(w)$,
$$
\Phi(Z\ox s)(X)=[Z,X]\cdot s
=\sum_{\mu}(w-{c}_{\la\mu})\pi_{\la\mu}
(Z\ox s)(X),
$$
where
$
{c}_{\la\mu}=-\frac{1}{2}[(\mu,\mu+2\rho)-(\la,\la+2\rho)
-(\al,\al+2\rho)].
$
\endproclaim

\demo{Proof}
Let $\{\xi_{\nu}\},$ resp. $\{\et_{\nu}\}$ be dual bases of $\gog_{-1}$,
resp\. $\gog_1$. The invariance of
the scalar product implies
$$
[Z,\xi_{\nu}]=\sum_a(Y_a,[Z,\xi_{\nu}]){Y}_a=\sum_a([Y_a,Z],\xi_{\nu})
{Y}_a
$$
$$
\Phi(Z\ox s)=
\sum_\nu\et_{\nu}\ox [Z,\xi_{\nu}]\cdot s
=\sum_\nu\et_{\nu}\ox\biggl(\sum_a([Y_a,Z],\xi_{\nu}){Y}_a\bigg)
\cdot s
=\sum_a[Y_a,Z]\ox {Y}_a\cdot s.
$$
Since  $Y_0=E$, the first term in the sum is
$
[Y_0,Z]\ox Y_0\cdot s=w\,Z\ox s
$
and for the rest we can use the definition of the Casimir
operator and its computation by means of highest weights,
together with
$$
\sum_{a'}Y_{a'}\,{Y}_{a'}\cdot(Z\ox s)=
\sum_{a'}(Y_{a'}{Y}_{a'}\cdot Z)\ox s+
\sum_{a'}Z\ox (Y_{a'}{Y}_{a'}\cdot s)+
2\sum_{a'}(Y_{a'}\cdot Z)\ox({Y}_{a'}\cdot  s)
$$
(notice $\cdot$  means the actions on
different modules used in the formula)
\qed\enddemo

\subhead\nmb.{4.4} Example\endsubhead
Let us compute now a simple case of the formula above which
will be needed below.  The special double commutator
terms  $[[X,\ta],\ta]$ from \nmb!{2.5} are appearing often in the
algorithm  mentioned in \nmb!{2.6}. We want to decompose
them into irreducible pieces.

Again, let $\al$ be the highest weight of $\gog_1$ considered
as $\gog_0^s$-module. By our conventions, it has the conformal weight
$1$. The tensor square
$\gog_1\ox\gog_1$ decomposes always into
symmetric and antisymmetric parts. But the symmetric square
decomposes in all but one cases
 into two components (an exceptional case being
projective structures, where is does not decompose). For our
purposes, it is sufficient to know that there is always a piece
in the decomposition with the highest weight $2\al$
(the Cartan product of $\gog_1$  with itself), denoted by
$\gog_1\x\gog_1$.

\proclaim{Lemma}
Let
$
\gog_1\ox\gog_1
=\oplus_{i=1}^3\VV_{\al_i}$
be the decomposition into irreducible components
with
$\VV_{\al_1}\simeq\x^2(\gog_1)$
and
$\VV_{\al_3}\simeq\La^2(\gog_1)$
(\,$\VV_{\al_2}$ is trivial in the projective case). Hence
$\al_1=2\al$.
Then there  exist real numbers $A_i,i=1,2,3,$
such that
$$
-{\frac{1}{2}}[[X,\ta],\ta](Y)=
\sum_{i=1}^3 A_i\pi_i[\ta\ox\ta](X,Y),
$$
where $X,Y\in\gog_{-1}$; $\ta\in\gog_1$,
and $\pi_i$ is the projection onto $\VV_{\al_i}$.
For $A_1$, we have $A_1={\frac{1}{2}}(|\al|^2+1)$.
\endproclaim

\demo{Proof}
This is the case $\VV_{\la}=\gog_1$ of lemma 4.3, so the
numbers $A_i$ are given by
$$
A_i=-{\frac{1}{2}}[{c}_{\al\al_i}-1],\quad i=1,2,3.
$$
In particular,
$
c_{\al,2\al}=-{\frac{1}{2}}[(2\al,2\al+2\rho)-2(\al,\al+2\rho)]=
-|\al|^2
$.
\qed\enddemo

In computations below, we shall use often the
constant $A_1$ but we shall see that its actual value
does not influence the explicit formula for standard operators,
because the constant $A_1$ will be absorbed by a renormalization of
the deformation tensor $\Ga$.

\head\nmb0{5}. $P$-module homomorphisms\endhead

Let us suppose, as in the previous section, that
$\gog$ is a complex $|1|$-graded
Lie algebra, $\goth p=\gog_0\oplus\gog_1$ and $\VV$ is
a (complex) irreducible $\goth p$-module.
The algebra $\gog_0$ splits into the sum of the commutative
$1$-dimensional ideal $\goth a$ and the semisimple
part $\gog_0^s$.

Using results obtained in the last two sections,
it is possible to construct a broad class of $\goth p$-homomorphisms
$\Phi$ from $\bar{J}^k\VV$ to a $P$-module $\VV'$,
where  $\VV'$
is a suitable irreducible component of the
$\gog_0$-module $\ox^k(\gog_1)\ox\VV$.
Let us recall that
there is a unique grading element $E\in\goth a$ for $\gog $ and
an invariant scalar product $(.,.)$ on $\gog$ is normalized
by the condition $(E,E)=1$.

Before stating the corresponding result, we shall prove a simple auxiliary
Lemma. A surprising and important fact coming from it
is the independence of the constants  $c_{j+1}-c_j$
of the chosen representations.

\proclaim{\nmb.{5.1} Lemma}
Let $\al$ be the highest weight of the $\gog_0^s$-module
$\gog_1$ and $\th$ one of its extremal weights.
For any  weight $\la$, let us define weights
$
\la_j=\la+j\th$, $j\in\NN
$,
and numbers
$$
c_j=c_{\la_{j}\la_{j+1}}=
-{\frac{1}{2}}\bigl[(\la_{j+1},\la_{j+1}+2\rho)-
(\la_{j},\la_{j}+2\rho)-(\al,\al+2\rho)\bigr].
$$
Then we have \roster
\item
$c_0=(\al,\rho)-(\th,\la+\rho);$
\item
$c_j-c_{j-1}=-|\al|^2;$
\item
$\sum_{j=0}^{k-1}c_j=k\,[(\al,\rho)-(\th,\la+\rho)-
{\frac{k-1}{2}}|\al|^2]$.
\endroster
\endproclaim

\demo{Proof}
By definition
$$
\align
c_0&=-{\frac{1}{2}}\bigl(
(\la+\th,\la+\th+2\rho)-(\la,\la+2\rho)-(\al,\al+2\rho)\bigr)=\\
&=(\al,\rho)-(\th,\la+\rho)-{\frac{1}{2}}(|\th|^2-|\al|^2).
\endalign
$$
The weight $\th$ lies in the $W$-orbit of $\al$, so they have
the same norm, and (1) follows.
Substituting $\la_j$ instead of $\la$, we get
$$
c_j=(\al,\rho)-(\th,\la+\rho)-j|\th|^2
$$
as well as the formula (2).
Using $c_j=c_0-j|\al|^2,$ we get
$$
\sum_{j=0}^{k-1}c_j=\sum_{j=0}^{k-1}(c_0-j|\al|^2)=
k\,c_0-{\frac{k(k-1)}{2}}|\al|^2.
\rlap{\qquad\qquad\qed}$$
\enddemo

\subhead \nmb.{5.2} The  algebraic criterion\endsubhead
We want now to prove that certain $G_0$-homo\-mor\-phisms are
in fact $P$-homomorphisms.
In \cite{CSS1}, the following algebraic condition for it
was proved, but in the case when the
invariant scalar product $(.,.)$ was equal to the Killing form
$B(.,.)$.
If the normalization of $(.,.)$ is different and if $\kappa$ is
a number such that $B(.,.)=\kappa(.,.)$, then it is easy to check
that all terms in the Lemma below are scaled
uniformly by the constant $\kappa ^k$, hence the condition
does not change.

\proclaim{Lemma}
Let $\VV$ and $\VV'$ be irreducible $P$-modules and
$\Ph\:\bar{J}^k (\VV)\to \VV'$ be a $\goth g_0$-module
homomorphism whose restriction to $\otimes^k
(\gog_{-1}^{*})\otimes \VV\subset \bar{J}^k(\VV)$ does not vanish.
Let us choose  any invariant scalar product $(.,.)$ on $\gog$
and let us use it to identify $\gog_1$ with $\gog_{-1}^*.$
Then $\Ph$ is a $P$-module homomorphism if and only if:
\roster
\item It factors through
the projection $\pi\:\bar{J}^k(\VV)\to
\otimes ^k
(\gog_{-1}^{*})\otimes \VV$;
\item $\Ph$ vanishes on the image of $\otimes ^{k-1}
(\gog_{-1}^{*})\otimes \VV$ in $\bar{J}^k(\VV)$
under the action of $\gog_1$, i.e\. for all
$Z,Y_1,\ldots,Y_{k-1}\in \goth g_1$, $v\in\VV$
we have
$$
\Phi\biggl(\sum_{i=0}^{k-1}
(\sum_\be Y_1\otimes\dots\otimes Y_i\otimes\et_\be\otimes
\bigl([Z,\xi_\be].(Y_{i+1}\otimes\dots\otimes Y_{k-1}\otimes v)
\bigr)\biggr)=0,
$$
\endroster
where $\et_\be$ and $\xi_\be$ are dual bases of $\goth g_1$ and $\goth g_{-1}$
with respect to the scalar product $(.,.)$
and the dot means the standard action of
an element in $\goth g_0$ on the argument.
\endproclaim

This criterion looks quite complicated.
Using results of Section \nmb!{4}, we can use it to prove easily the existence
of a broad class of $P$-modules homomorphisms.

\proclaim{\nmb.{5.3} Corollary}
Let $\VV_{\la}$ be an irreducible $\gog_0^s$-module and let $\al$ be the
highest weight of the irreducible $\gog_0^s$-module $\gog_1$.

Let us suppose that an extremal weight $\th$ of $\gog_1$ and an
positive integer $k$ is chosen in such a way that $\mu=\la+k\th$  is
dominant. Let $\pi:\ox^k\gog_1\ox\VV_{\la}\to\VV_{\mu}$
be the projection on the unique irreducible component of the
product with highest weight $\mu$ (see Theorem \nmb!{3.5}).

Then there is a unique value for the generalized conformal weight
$w$ such that $\pi$ defines a $P$-homomorphism from
$\bar{J}^k(\VV_{\la}(w))$ to $\VV_{\mu}(w+k).$ The value of
that conformal weight is given by
$$
w=(\al-\th,\rho)-{\frac{k-1}{2}}(|\al|^2+1)-(\th,\la),
$$
where $\rho$ is half the sum of positive roots for
$\gog_0^s$.
\endproclaim

\demo{Proof}
Let us first recall  the  construction of the projection  $\pi$.
If $\la_{k'}=\la+k'\th$, $k'=0,\ldots,k$,
the projections $\pi_{k'}$, $k'=1,\ldots,k$, are defined inductively as
the projections from $\gog_1\ox\VV_{\la_{k'-1}}$
onto the unique irreducible component $\VV_{\la_{k'}}$
with highest weight $\la_{k'}$. The projection $\pi$ is
given by the formula
$$
\pi(Z_{1}\otimes\dots\otimes Z_{k}\otimes v_{})=
\pi_k(Z_1\ox\pi_{k-1}(Z_2\ox\ldots\pi_1(Z_k\ox v)\ldots )),
$$
where  $Z_1,\ldots,Z_k\in\gog_1$, $v\in\VV_{\la}$.

To prove the theorem, we have
to verify that with the choice of the weight $w$
above, the condition in Lemma \nmb!{5.2} is satisfied.
So we want to find $w$ in such a way that
for all $Z,Z_1,\ldots Z_{k-1}\in\gog_1$, $v\in\VV_{\la}$,
$$
\pi\biggl(\sum_{i=0}^{k-1}
\sum_\be Z_1\otimes\dots\otimes Z_i\otimes\et_\be\otimes
\bigl([Z,\xi_\be].(Z_{i+1}\otimes\dots\otimes Z_{k-1}\otimes v_{})
\bigr)\biggr)=0,
$$
where $\et_\be$ and $\xi_\be$ are dual bases of
$\goth g_1$ and $\goth g_{-1}$
with respect to the  product $(.,.)$.
Let us recall the notation
$c_j=c_{\la_j,\la_{j+1}}$ from Lemma \nmb!{5.1}.

By Lemma \nmb!{4.3}, applied to elements from
$\VV_{\la_{k-1-i}}(w+k-1-i)$, we have
$$
\multline
\pi_{k-i}\biggl(
\sum_{\be}\et_\be\otimes\pi_{k-i-1}
\bigl([Z,\xi_\be].(Z_{i+1}\otimes\pi_{k-i-2}
(\ldots\otimes\pi_1( Z_{k-1}\otimes v_{})\ldots))
\bigr)\biggr)=\\
\pi_{k-i}\biggl(
\sum_{\be}\et_{\be}\otimes
\bigl([Z,\xi_\be].
(\pi_{k-i-1}
(Z_{i+1}\otimes\pi_{k-i-2}
(\ldots\otimes\pi_1( Z_{k-1}\otimes v_{})\ldots)))
\bigr)\biggr)=\\
\left(w+k-1-i-{c}_{k-1-i}\right)\pi_{k-i}
\bigl(Z\ox\pi_{k-i-1}
(Z_{i+1}
(\ldots\otimes\pi_1( Z_{k-1}\otimes v_{})\ldots))
\bigr).
\endmultline$$
Due to the fact that all images of $\pi^j$ belong to
$\odot^j\gog_1\ox\VV_{\la}$, $j=1,\dots,k$, all elements
$$
\pi(Z_1\otimes\dots \ox Z_i\otimes Z\ox Z_{i+1}\otimes\dots\otimes
Z_{k-1}\otimes v_{}));\;i=0,\ldots,k-1
$$
coincide. It is hence sufficient to find $w$ so that
$$
kw+{\frac{k(k-1)}{2}}-\sum_{j=0}^{k-1}
{c}_{k-1-j}=0.
$$
To get the value for $w$, it is sufficient
to use Lemma \nmb!{5.1} (note that $|\al|=|\th|$).
\qed\enddemo

\head\nmb0{6}. Standard operators\endhead
\subhead\nmb.{6.1} A construction of invariant operators\endsubhead
As described in Section \nmb!{2}, the $P$- module homo\-morphisms constructed
in the last Section define
invariant differential operators. We can now summarize
the whole construction and the data needed for it.
Let us  return to the situation of Section \nmb!{2}
with a given $|1|$-graded (real) simple Lie algebra $\gog,$
the corresponding groups $P\subset G,$ $G_0$,
and a principal fiber bundle $\G$ over $M$ with
a given Cartan connection $\om$.

The complexification $\gog^{\CC}$ is a complex semisimple
$|1|$-graded Lie algebra and
$\gog_j=\gog\cap\gog^{\CC}_j;\;j=-1,0,1$.
Any (complex) irreducible $P$-module $\VV$ is an irreducible
$\gog_0$-module
as well as $\gog_0^{\CC}$-module.
They are characterized by an integral dominant weight for
$(\gog_0^s)^{\CC}$ and the (generalized) conformal weight $w$.
The tensor product
$\gog_1\ox_{\RR}\VV$ is isomorphic to
$\gog_1^{\CC}\ox_{\CC}\VV$, the same is true for iterated
tensor products.
The  space $\gog_1^{\CC}$ is an irreducible module for
$\gog_0^s$ with a highest weight $\alpha$.

Suppose that we have chosen the following data:
An irreducible module $\VV_{\la}$ for $\gog_0^s$,
a 'direction' $\th$, which is an extremal weight of
the $\gog_0^s$-module $\gog_1^{\CC}$, and
a positive integer $k$,  such that
$\mu=\la+k\th\in\La^+$.

Let $\pi$ be the projection to the unique irreducible component
of the $\gog_0^s$-module
$\ox^k\gog_1\ox\VV_{\la}$ with the highest weight $\mu=\la+k\th$ (cf\.
Theorem \nmb!{3.5}), and let $w$ be the corresponding (generalized) conformal
weight from Corollary \nmb!{5.3}.
Then the operator
$$
D\equiv D(\la,\th,k)=\pi\o(\nabla^{\om})^k:
{\Cal C}^{\infty}(P,V_{\la}(w))^{P}
\>{\Cal C}^{\infty}(P,V_{\mu}(w+k))^{P},
$$
is an invariant differential operator of order $k$.

\subhead\nmb.{6.2} Standard operators\endsubhead
We have defined above a certain class of operators which were
proved to be invariant. There is a traditional division of
invariant operators into two classes --- standard and nonstandard
ones. We would like to show now that the operators constructed
above include almost the whole class  of
so called  standard operators.

(Fundamental)
standard operators were originally defined in the homogeneous situation
(on generalized flag manifolds $G/P,$ with $G$ complex simple
and $P$ parabolic). In the Borel case, the classification
of all invariant differential operators was given (in the dual language
of homomorphism between Verma modules) by Bernstein, Gelfand and
Gelfand, see \cite{BGG}. They are all defined uniquely by their
source and target (up to a constant multiple) and they
are precisely all operators forming the so called BGG resolutions.
For a general parabolic, the BGG resolutions are also well known
but the class of invariant operators corresponding to individual
arrows in them --- they are called
(fundamental) standard operators ---
is no more the complete set of invariant operators.
There exist also the so called non-standard operators.
To show a relation of our invariant operators $D(\la, \th, k)$ to the
standard operators, we need just their following simple
property (more details can be found e.g\. in \cite{BasE}, \cite{Go}).

Suppose that a Cartan subalgebra $\goth h$ in $\gog^{\CC}$
and the set of simple roots is chosen
in such a way that $E\in\goth h$ and that all positive
spaces are contained in $\gog_0^{\CC}\cap\gog_1^{\CC}.$
Irreducible representations of $\gog_0^{\CC}$ can be characterized
by their highest weight, considered as an element in $\goth h^*,$
such that its restriction to $(\goth h)^s=\goth h\cap (\gog_0^{\CC})^s$
is dominant. This carries information both on the highest weight
for $(\gog_0^{\CC})^s$ and on a generalized conformal weight.
For any such $\La\in\goth h^*,$ the symbol $V_{\La}$ denotes
a homogeneous bundle
given by the irreducible representation of $\gog_0^{\CC},$ corresponding
to this highest weight.
The Weyl group $W$ of $\gog^{\CC}$ has a structure of a directed
graph which is directly related to existence of invariant
operators.

The property we need is the following. If
$D:\Ga(V_{\La})\to\Ga(V_{\La'})$ is a standard invariant operator,
then there is a positive root $\Theta$ for $\gog^{\CC}$
such that $\si_{\Theta}(\La+\De)=\La'+\De$, where
$\si_{\Theta}$ is the reflection with respect to $\Theta$ and
$\De$ is a half-sum of positive roots for $\gog^{\CC}$.
Consequently, we have also $|\La+\De|^2=|\La'+\De|^2$.
Before going further, we need two simple auxiliary lemmas.

\proclaim{\nmb.{6.3} Lemma}
Let $\gog$ be a complex $|1|$-graded Lie algebra,
$S=\{\al_i\}_{i=0}^m$ the set of its
simple roots with its numbering chosen in such a way that $\al_0$
is the crossed simple root.
Let $\{\pi_i\}$ be the corresponding set of fundamental weights.

Then we have
\roster
\item If $\La$ is the highest weight of an irreducible $\gog_0$-module $\VV$,
then its conformal weight is equal to $w=\La(E)$.
\item
The root space $\gog_{\al}$ belongs to $\gog_j$, $j=-1,0,1$,
if and only if $a_0=j$, where $a_i$ are coefficients  in the decomposition
$\al=\sum_{i=0}^ma_i\al_i$.
\item For any weight $\La\in\goth h^*$,
we have
$
(\pi_0,\La)={\frac{|\al_0|^2}{2}}\La(E),
$
where $E$ is the grading element.
\item
Let us consider two weights $\La$, $\La'$ and a number $a$
such that
$|\La|^2=|\La'|^2$, $|\La+a\pi_0|^2=|\La'+a\pi_0|^2 $
and $(\La-\La',\pi_0)\not=0$.
Then   $a=0$.
\endroster
\endproclaim

\demo{Proof}
(1)
If $v$ is a highest weight vector for $\VV,$ then
$E\cdot v=\La(E)v,$ but by definition $E\cdot v=w\, v.$

(2)
This is a special case of a simple general statement valid
for all $|k|$-graded Lie algebras. The reason is that all simple roots but
$\al_0$ are in $\gog_0$, while $\al_0$ generates $\gog_1$.

(3)
There is an element $H\in\goth h$
such that $(\pi_0,\La)=\La(H)$
for all $\La\in\goth h^*$.
Then for all $j=1,\ldots,m$, we have $0=(\pi_0,\al_j^\vee )=\al_j^\vee
(H)$,
where $\al_j^\vee ={\frac{2\al_j}{|\al_j|^2}}$.
The element $H$ is orthogonal to all roots of $\gog_0$, hence
it is a multiple of $E$ (which has the same property).
To check the multiple, it is sufficient to note that $\al_0(E)=1,$
because the conformal weight for $\gog_1$ is $1.$
\item{4)}
The last property follows from
$$
|\La+a\pi_0|^2-|\La'+a\pi_0|^2=
2a(\La-\La',\pi_0).\rlap{\qquad\qquad\qquad\qed}
$$
\enddemo

As a consequence, we get the following interesting fact.

\proclaim{\nmb.{6.4} Lemma}
In the setting of \nmb!{6.1}, let $\la,\la'$ be two dominant integral
weights for $\gog_0^s.$
Suppose that there are two nontrivial standard invariant differential
operators $D,\tilde{D}$ of order $k>0$ such that
$$
D:\Ga(\VV_{\la}(w))\to\Ga(\VV_{\la'}(w+k));\;
\tilde{D}:\Ga(\VV_{\la}(\tilde{w}))\to\Ga(\VV_{\la'}
({\tilde{w}}+k)).
$$
Then $w=\tilde{w}.$
\endproclaim

\demo{Proof}
Let  $\La$, $\La'$, $\tilde{\La}$, $\tilde{\La'}$ be in turn highest
weights from $\goh^*$ for irreducible representations
$$
\VV_{\la}(w),\VV_{\la'}(w+k),\VV_{\la}(\tilde{w}),
\VV_{\la'}(\tilde{w}+k).
$$
If $\De$ is the half-sum of positive roots for $\gog,$ then
existence of $D,\tilde{D}$ implies (see \nmb!{6.2}) that
$$
|\La+\De|^2=|\La'+\De|^2;\;
|\tilde{\La}+\De|^2=|{\tilde{\La}}'+\De|^2.
$$
The differences $\tilde{\La}-\La$, $\tilde{\La'}-\La'$
annihilate $\goh^s$, hence there are numbers $a$, $a'$ such that
$\tilde{\La}-\La=a\pi_0$; $\tilde{\La'}-\La'=a'\pi_0$.
But
$$
a\pi_0(E)=(\tilde{\La}-\La)(E)=\tilde{w}-w=  (\tilde{\La'}-\La')(E)=a'
\pi_0(E),
$$
hence $a=a'$.
Moreover,
$(\La-\La')(E)=k>0$, hence $(\La-\La',\pi_0)\not=0$. Now,
Lemma \nmb!{6.3} implies that $a=0$.
\hfill\qed\enddemo

\proclaim{\nmb.{6.5} Theorem}
Let $D$ be a standard invariant differential operator
acting between sections of $V_{\La}$ and $V_{\ti\La}.$
Let $\Theta\in\goh^*$ be a positive root of $\gog$ such that
$\ti\La+\De=\si_{\Theta}(\La+\De).$
Denote by $\th$ the restriction of $\Theta$ to $\goh^s$
and by $\la$ the restriction of $\La.$

Then $\th$ is a weight of $\gog_0^s$-module $\gog_1$ and
the number
$k=2(\La+\De,\Theta)/(\Theta,\Theta)$
is a positive integer.

If moreover the weight  $\th$ is an extremal weight of $\gog_1,$
then the operator $D(\la,\th,k)$ defined in 6.1 coincides
(up to a multiple) with the operator $D$ on sections of 
the homogeneous bundle $V_{\La}$.
\endproclaim

\demo{Proof}
The root $\Theta$ is a positive root of $\gog.$ Consequently,
the value of $\Theta(E)$ is either 0 or 1.
By the properties of standard operators (see \nmb!{6.2}), we have
$$
\ti\La-\La = k\Theta,
$$
where $k=2(\La+\De,\Theta)/(\Theta,\Theta)$ must be an integer. Because any
differential operator must increase (generalized) conformal
weight (which is given by evaluation of the highest weight on $E$),
 the value $\Theta(E)$ cannot vanish. Hence $\Theta(E)=1$
and $k>0.$

If we denote by $\la$, resp\. $\ti\la$, the restrictions of $\La$,
resp\. $\ti\La$ to $\goh^s$, then
we have also the relation
$$
\ti{\la}=\la+k\th.
$$
Hence the operators $D$ and $D(\la,\th,k)$ act between  the same
$\gog_0^s$ bundles  and
they are both invariant. By Lemma \nmb!{6.4}, their conformal weights
coincide as well. Now, the standard operators are completely
defined by their domains and targets up to multiples, see \cite{BC}, and $D$
and $D(\la,\th,k)$ differ at most by a constant
multiple.\hfill\qed
\enddemo

\subhead\nmb.{6.6} Remark
\endsubhead
We have just seen that our construction gives all standard
invariant operators for those AHS structures, for which
the set of weights of $\gog_1^{\CC}$ is just one orbit
of the Weyl group. This is true for all cases with two exceptions
--- the odd dimensional conformal case and the symplectic case.

There is indeed an exceptional set of standard operators for AHS
structures which do not have a simple description
of the form $D(\lambda,\theta,k)$ constructed above.
A typical  example is the case of odd conformal structures
and the operators in the middle of the BGG resolution.
These are operators acting between sections
$\Gamma({\Bbb V}_{\lambda}(w))$ and
$\Gamma({\Bbb V}_{\lambda}(w'))$.
The representation ${\Bbb V}_{\lambda}$ of the semi-simple part of $G_0$
is the same for
the source and the target, they differ only by their conformal
weights. They correspond to the case of operators $(\lambda,\theta,k)$,
where $\theta$ is the zero weight of $\goth g_1$. In this case, however,
the isotypic component ${\Bbb V}_{\lambda}$ appears in
$\otimes^k({\goth g}_1)\otimes{\Bbb V}_{\lambda}$ with higher
multiplicities.

In general, the BGG sequence of a representation ${\Bbb V}$ of
${\goth g}$ can be realized using the twisted (${\Bbb V}$-valued)
de Rham sequence.
In the particular case of the BGG sequence of the basic spinor
representation ${\Bbb S}$ of $\goth g=Spin(2n+2,\CC)$,
the middle operator corresponds to a second order operator $D$ between
$\Gamma({\Bbb V}_{\lambda}(n-1/2))$,
and $\Gamma({\Bbb V}_{\lambda}(n+3/2))$,
where $\lambda=(3/2,\ldots,3/2)$.
There are 3 pieces in the decomposition of the tensor
product $\otimes^2({\goth g}_1)\otimes\Gamma({\Bbb V}_{\lambda})$,
corresponding to sequences of weights $(\lambda,\sigma,\lambda)$
with $\sigma_1=(5/2,3/2,\ldots,3/2)$; $\sigma_2=(3/2,\ldots,3/2)$;
$\sigma_3=(3/2,\ldots,3/2,1/2)$.
It can be shown by methods described in \cite{CSS4}, \cite{B}, (see also
\cite{Sev}) that the corresponding standard operator is given by
$\pi\circ(\nabla^{\gamma})^2$,
where the projection $\pi$ is equal to
$\pi=\pi_2+1/4\pi_3$, where $\pi_j$ are defined as projections
to irreducible pieces corresponding to the sequences with $\sigma_j$.
The form of the operator $D$ is hence more complicated,
it has the form
$$
D\, t=\pi_2[(\nabla^{\gamma})^2 t-(1/2)\Gamma\otimes t]
+1/4\pi_3[(\nabla^{\gamma})^2 t-2\Gamma\otimes t].
$$
So it is clear that its formula has no more the simple universal form
$D\, t=\pi((\nabla^{\gamma})^2 t+\Gamma\otimes t])$ of the second order
standard operators deduced below, see \nmb!{7.11}.

\head\nmb0{7}. Explicit formulae for standard operators\endhead

\subhead\nmb.{7.1} Obstruction and correction terms\endsubhead
An algorithm for computation of $(\na^{\om})^k$ in terms
of the principal connection $\na^{\ga}$ and
its deformation tensor $\Ga$ was given in \cite{CSS1}, Sec\. 4.
The formulae for obstruction terms (important for existence proofs)
as well as for correction terms (important for explicit
description of operators) become quickly very complicated.
Using explicit description of the homomorphism $\Phi$ in
Section \nmb!{4} by means of Casimir operators, it is possible to
simplify the algorithm substantially and to
get quite explicit formulae for the coefficients in general
correction terms for the invariant operators constructed in the previous
section.
It is quite remarkable that  coefficients in the final formula
for curvature correction terms
do not depend on a choice of
a representation $V_{\la}$
as well as on a choice of a particular AHS
structure! They depend only on the order of the operator.

Let us first simplify the algorithm given in \cite{CSS1}.
Let $k$ be a fixed integer and let us consider an operator
$D=\pi\o(\nabla^{\om})^k$,
where the projection $\pi$ of
$\ox^k(\gog_1^{\Bbb C})\ox V_{\la}$ onto one of its irreducible
components is determined by a chain
of dominant weights, as described in Section \nmb!{3}.
Knowing  highest weights of all intermediate irreducible
components in the chain of projections,
Lemma \nmb!{4.3} can be used to compute the values of the homomorphism  
$\Phi$ on all terms in the algorithm. The same is true for
the action of the double commutator term $[[X,\ta],\ta]$
(see Example \nmb!{4.4}).
This makes it possible to evaluate, in principle, all terms in the
expansion. But the result is still quite complicated.

A considerable simplification in the algorithm can be achieved,
if we restrict ourselves to the symmetric case, i.e\. if
the image of $\pi$ is a subspace of
$\odot^k(\gog_1^{\Bbb C})\ox V_{\la}$.
Then many multiple tensor products contained in various terms
of the formula may be reordered and combined together.
Any term of the formula is then just a symmetric tensor product
of a power of $\ta$, suitable powers of $\Ga,$ its
covariant derivatives and a covariant derivative of the section
$s.$ A problem to be solved is whether there is a way how to compute
effectively coefficients in the corresponding linear combination
of such terms.

An additional simplification can be achieved in the case, when
we know which summand in the description of the action of the
double  commutator (Lemma \nmb!{4.4}) is really appearing in various
terms. Such information is available in the case of the operators
$D(\la,\th,k)$ constructed above.
In this case, we may use properties of the decomposition of the
tensor product
$\ox^k(\gog_1^{\Bbb C})\ox V_{\la}$ proved in Section \nmb!{3}
to get an explicit form of the operator. Before tackling the main Theorems
\nmb!{7.4} and \nmb!{7.9}, we discuss the low order cases.

\subhead\nmb.{7.2} The first  order operators\endsubhead
Using results from \cite{CSS1}, see \nmb!{2.4}, and Lemma \nmb!{4.3},
 we get immediately the existence and an explicit
form of the 1st order operators.

\proclaim{Corollary}
Let $\VV_{\la}$ be an irreducible representation of
$(\gog_0^{\CC})^s$ and
$\VV_{\mu}$ be an irreducible component of the
product
$\gog_1\ox \VV_{\la}$. Let $\pi=\pi_{\la\mu}$
be the corresponding projection. Then
$$
\pi(\na^{\om}(p^*t))=\pi[p^*(\na^{\ga} t)+(c_0-w)\ta\ox t]
$$
where $c_0=c_{\la\mu}$ are the constants from \nmb!{4.3}.

In particular, there is the unique value $w=c_0$ of the conformal weight
for which the projection defines a first order invariant
operator $D\, t=\pi[p^*(\na^{\ga}) t]$.
\endproclaim

Operators of this type were introduced in conformal case
in paper \cite{SW} and are now standardly called
generalized gradients or Stein--Weiss operators (see e.g\.
\cite{Bra}).
The result above was proved in the conformal case
by Fegan (see \cite{F}).
He gave the first systematic classification of such operators.
The theorem above treats completely all first order operators
for all AHS structures (note that in odd conformal case,
the class of them includes also certain exceptional
standard operators of first order not covered by the class
of operators $D(\la,\th,k)$, e.g\. the one in the middle in
the de Rham resolution).

\subhead\nmb.{7.3} The second  order operators\endsubhead
In a similar way, we can use the first order formula, the algorithm leading
in \cite{CSS1} to the formula in \nmb!{2.5}, and Lemma \nmb!{4.2}, in order
to compute explicitly the form of the second order invariant
differential projected to an irreducible component
given by a sequence of dominant weights
$\underline{\la}=(\la_0,\la_1,\la_2)$. Let $\pi$ be the corresponding
projection.

\proclaim{Corollary}
Using notation of Example \nmb!{4.4} and Lemma \nmb!{5.1}, we have
$$
\align
\pi\bigl[\bigl( (\nabla^\om)^2 (p^*t) \bigr)]=\ &
\pi[p^*( (\nabla^{\ga})^2t)+
(c_0-w)\Ga\ox p^*t+\\
&(c_0-w)\ta\ox p^*(\nabla^{\ga}t)+
(c_1-w-1)p^*(\nabla^{\ga}t)\ox\ta+\\
&(c_0-w)(c_1-w-1)\ta\ox\ta\ox t -
\sum_{i=1}^3A_i\pi_i(\ta\ox\ta\ox t)\bigr].
\endalign
$$
\endproclaim

The most complicated term to compute is clearly the last one
coming from the double commutator term. To understand that term, one
has to understand well the relation among the chosen projection $\pi$
defined by the chain of weights $\underline\la$ and the projections
$\pi_i$ coming from the splitting $\gog_1\ox\gog_1$ into symmetric
and antisymmetric parts. We shall see that for operators
$D(\la,\th,k)$, this relation can be understood and the formula
above can be simplified further.

The operators $D(\la,\theta,2)$ are invariant for a unique value
for the (generalized) conformal weight, cf\. \nmb!{6.1}.
It is immediate to check that it is just given
by the requirement that the sum of coefficients at  terms linear in $\ta$
vanishes. It is also possible to verify directly that then the coefficient
at the term of second order in $\ta$ vanishes as well.

We shall now follow line of reasoning suggested in \nmb!{7.1} and
we shall develop an effective procedure for explicit
description of all operators $D(\la,\th,k)$.

\proclaim{\nmb.{7.4} Theorem}
Let $A_1$ be the number defined in Example \nmb!{4.4}.
The value of the operator
$D(\la,\theta,k)t(u)=\pi_k\circ((\na^{\om})^k(p^*t))(u)$
constructed in \nmb!{6.1} expands into a sum of the form
$$
\sum
a^{k,j}_{s_0,\ldots,s_m}
\pi_{k}[\ta^j\odot\Ga^{s_0}\odot (\na\Ga)^{s_1}\odot\ldots
\odot(\na^m\Ga)^{s_m}\odot \na^i\,t](u),
$$
where the summation goes over
$$
j,s_i\in\{0,1,2,\ldots\}\text{\ such that\ }j+\sum_{i'=0}^m
s_{i'}(i'+2)+i=k, 
$$
$a^{k,j}_{s_0,\ldots,s_m}\in \RR$, $\ta(u)\in\gog_1^{\Bbb C}$, and
$$\gather
\ta^j=\odot^j\ta,\quad [\na^i\,t](X_1,\ldots,X_i)=
p^*\na^{\ga}_{X_i}\ldots\na^{\ga}_{X_1}t,
\\
[\na^\ell\Ga](X,Y,X_1,\ldots,X_\ell)=
[p^*\circ\na^{\ga}_{X_\ell}\ldots\na^{\ga}_{X_1}(\Ga)]
(X,Y).
\endgather
$$
The expressions
$$
F^{k} t(u) := \pi_{k}[(\nabla^\om)^{k} (p^*t)](u)\in
\odot^{k}(\gog_1^{\Bbb C})\ox V_{\la}
$$
are given by  recursive formulae
$$\align
F^0t(u) &= p^*t(u)\\
F^{k+1}t(u) &= [S_{\la+\ta}](F^{k}t(u)) +
 [S_{\nabla}](F^{k}t(u))+
 [S_{\Ga}](F^{k}t(u)).
\endalign
$$
The individual transformations
$S_{\la+\ta},S_\nabla$ and $S_\Ga$ act as follows:
$$\align
&S_{\la+\ta}[\pi_k(\ta^{j-1}\odot\om_{k-j+1})]=
({c}_{k}-k+(j-1)A_1-w)
\pi_{k+1}(\ta^{j}\odot\om_{k-j+1}\ox t)        ;
\\
&\qquad\text{where $\om_{k-j+1}\in\odot^{k-j+1}(\gog_1^{\Bbb C})\ox V_{\la};\;
{c}_{k}={ c}_{\la_{k},\la_{k+1}};\;\la_k=\la+k\theta,\,j>1$.}
\\
&\aligned
S_{\na}[\pi_{{k}}
&(\ta^j\odot\Ga^{s_0}\odot (\na\Ga)^{s_1}\odot\ldots
\odot(\na^m\Ga)^{s_m}\odot \na^i\,t)]=\\
=\
&s_0[\pi_{{k+1}}(\ta^j\odot\Ga^{s_0-1}\odot (\na\Ga)^{s_1+1}\odot\ldots
\odot(\na^m\Ga)^{s_m}\odot \na^i\,t)]+\\
&+\ldots+\\
&s_m[\pi_{{k+1}}(\ta^j\odot\Ga^{s_0}\odot\ldots
\odot(\na^m\Ga)^{s_m-1}\odot
(\na^{m+1}\Ga)\ox \na^i\,t]+\\
&[\pi_{{k+1}}(\ta^j\odot\Ga^{s_0}\odot (\na\Ga)^{s_1}\odot\ldots
\odot(\na^m\Ga)^{s_m}\odot \na^{i+1}\,t)].
\endaligned
\\
&S_{\Ga}[\pi_{{k}}(\ta^{j+1}\odot\om_{k-j-1})]=
(j+1)\pi_{{k+1}}(\ta^{j}\odot\Ga\odot\om_{k-j-1})        ;
\\
&\qquad\text{where $\om_{k-j-1}\in\odot^{k-j-1}(\gog_1^{\Bbb C})\ox
V_{\la}$.}
\endalign$$
\endproclaim

\demo{Proof}
In \cite{CSS1, 4.9}, we have described an algorithm to inductively
compute the difference
$(\nabla^\om)^k (p^*t)- p^*((\nabla^\ga)^k t)$
as a sum of correction and obstruction terms.
Computing instead of that difference the value of
$F^kt(u):=(\nabla^\om)^k (p^*t)$ inductively, the results of
\cite{CSS1, 4.9} read as follows: The expression $F^kt(u)$, evaluated
at $k$ arguments from $\gog_{-1}$, expands into a sum of terms of the
form 
$$
a\la^{(t_1)}(\be_1)\dots\la^{(t_i)}(\be_i)p^*(\nabla^\ga)^jt
$$
where $a$ is a scalar coefficient,
the $\be_\ell$ are iterated brackets involving
some arguments $X_\ell\in\gog_{-1}$, the iterated covariant
differentials $(\nabla^{\ga})^r\Ga$ evaluated
on some $X$'s, and $\ta$'s. Exactly the first $t_j$ arguments
$\row X1{t_j}$ are evaluated after the action of
$\la^{(t_j)}(\be_j)$, the other ones appearing on the right are
evaluated  before.
For $k=1$, we have
$$
F^1t(u)(X_1) = p^*((\nabla^\ga) t)(u)(X_1) + [X_1,\ta](p^*t)(u).
$$
Inductively,
$$\align
F^kt(u)(\row X1k) =\ & \bigl(\la^{(k-1)}([X_k,\ta(u)])
F^{k-1}t(u)\bigr)(\row
X1{k-1}) + \\
&
\tilde S_{\ta}(F^{k-1}t(u))(\row X1{k}) +\\
&
\tilde S_{\nabla}(F^{k-1}t(u))(\row
X1{k})+\\
&
\tilde S_{\Ga}(F^{k-1}t(u))(\row X1{k}) .\\
\endalign$$
where $\la^{(k-1)}$ is the obvious tensor product representation on
$\ox^{k-1}\gog\ox V_{\la}$ and the individual transformations
$\tilde S_\ta$, $\tilde S_\nabla$, and $\tilde S_\Ga$ act as follows.
\roster
\item The action of $\tilde S_\ta$ replaces each summand
$a\la^{(t_1)}(\be_1)\dots\la^{(t_i)}(\be_i)p^*(\nabla^\ga)^jt$ by
a sum with just one term for each occurrence of $\ta$ where this
$\ta$ is replaced by $[\ta,[\ta,X_k]]$ and the coefficient $a$ is
multiplied by $-1/2$.
\item $\tilde S_\nabla$  replaces each summand in $F^{k-1}$ by a sum
with just one term for each occurrence of $\Ga$ and its
differentials, where these arguments are
replaced by their covariant derivatives $\nabla^\ga_{X_k}$, and
with one additional term where $(\nabla^\ga)^jt$ is replaced by
$\nabla^\ga_{X_k}((\nabla^\ga)^jt)$.
\item $\tilde S_\Ga$ replaces each summand by a sum with just one term
for each occurrence of $\ta$ where this $\ta$ is replaced by
$\Ga(u).X_k$.
\endroster

Now we are going to specialize these results to the case we are
interested in here: 
Under the assumptions of the theorem, which we want to prove,
the image of the projection $\pi$ is included in
$\odot^{k}(\gog_1^{\Bbb C})\ox V_{\la}$ hence
order of factors in the multiple tensor product does not matter.
Consequently all $\ta$'s can be shifted to the front of the product,
derivatives of $\Ga$ can be reordered as indicated above, and
all derivatives of $t$ can be put to the end of the expression.
Terms $\na^l\Ga$ can be hence interpreted as elements of
$\odot^{l+2}(\gog_1^{\Bbb C})\ox V_{\la} $ and
$\na^i\,t$ can be substituted by its symmetrization in
$\odot^{i}(\gog_1^{\Bbb C})\ox V_{\la}$.
We have already seen that the expression $F^1t$ has the required form
(see \nmb!{7.2}).
Using Casimir operators, we can now express the algorithm described
above in the following way.

Suppose (by induction) that the term $F^{k}$ has already been written
in the form given in the theorem.  The action of an element
$[X_{k+1},\ta(u)]$
on $F^{k}t(u)$
can be computed by Lemma \nmb!{4.3},
because we know that  $F^kt(u)$ belongs to the image of $\pi_{k}$,
which is, by assumption, an irreducible representation with
the highest weight $\la_{k}$. The result is
$(c_k-w-k)F^kt(u)$.

%
The action of $\tilde S_{\ta}$  was a replacement of
$\ta$ at all $j-1$ places in the expression by
$-1/2[\ta,[\ta,X_k]]$
Applying the projection $\pi$ and using the result of Example \nmb!{4.4} and
\nmb!{3.6},
only the first  part in the decomposition
of $\ta\ox\ta$  survives and the result is
the same term containing one more $\ta$ multiplied by $(j-1)A_1$.
Adding both contributions, we get the action of  $S_{\la+\ta}$.

The action of $\tilde S_{\na}$  is just a derivation and action
of $\tilde S_{\Ga}$ is a substitution of $\Ga$ instead of $\ta$, so we
arrive directly at the description of $S_\na$ and $S_\Ga$ in the
theorem. 

The fact that $F^k$ has the required form follows
from the above description of the operators $S_{\la+\ta}$, $S_\na$, $S_\Ga$
by induction.
\hfill\qed
\enddemo

Looking at the action of the individual transformations and at the form of the
expansion, we get immediately the following algorithm for
the unknown coefficients.

\proclaim{\nmb.{7.5} An algorithm for expansion coefficients}
The coefficients $a^{k+1,j}_{s_0,\ldots,s_m}$ in theorem \nmb!{7.4}
satisfy the following recursive relations.
$$
\align
a^{k+1,j}_{s_0,\ldots,s_m}=\ &
(1-\de_{j,0})a^{k,j-1}_{s_0,\ldots,s_m}
({c}_k-k+(j-1)A_1-w)\\
&+a^{k,j}_{s_0,\ldots,s_m}\\
&+(1-\de_{s_0,0})(j+1)a^{k,j+1}_{s_0-1,s_1,\ldots,s_m}+\\
&+(1-\delta_{s_1,0})(s_0+1)a^{k,j}_{s_0+1,s_1-1,\ldots,s_m}+\\
&+\ldots+\\
&+(1-\delta_{s_{m},0})(s_{m-1}+1)
a^{k,j}_{s_0,\ldots,s_{m-2},s_{m-1}+1,s_m-1}.
\endalign$$
\endproclaim

\subhead\nmb.{7.6} Constants $\tilde c_k$\endsubhead
In the algorithm above, the value $c_k-k+jA_1-w$
has frequently appeared.
It will be convenient to change the definition of constants
$c_j$ and to
define new shifted constants $\ti{c}_j$ instead.
Let us define them by
$$
\tilde{c}_j=c_0-j\,A_1.
$$
Then $c_k-k+j\,A_1-w=c_0-kA_1-(k-j)A_1-w=\tilde{c}_k-(k-j)A_1-w$.

Note for future use that the differences
$\tilde{c}_j-\tilde{c}_k=(k-j)A_1$ are  always multiples
of $A_1$.

\subhead\nmb.{7.7} Constants $B^m_{(s_0,\ldots,s_m)}$\endsubhead
As the last item in the preparation of an explicit computation
of the coefficients in the expansion, we are going to define inductively
the following parametric system of  constants
$B^n_s,$
where $n\geq 0$ is an integer,
$s=(s_0,s_1,s_2,\ldots)$ is a sequence of non-negative integers
with  a finite number of nonvanishing elements.
We shall often write $s=(s_0\ldots s_m)$ by cutting the sequence
at the last nontrivial entry; $(0)$ will denote the sequence
$(0,0,\ldots)$.
For any finite sequence of integers $s$ ,we shall use
two integers $|s|$, $[s]$ associated with $s$,  defined by
$$
|s|=\sum_0^{\infty}s_i\text{\ and\ }[s]=\sum_0^{\infty}s_i(i+1).
$$

Symbols $\si_i,\,i=0,1,\ldots$, will be used for special sequences
of integers defined by
$$
\si_0=(1,0,\ldots);\;
\si_1=(-1,1,0,\ldots);\;
\si_2=(0,-1,1,0,\ldots);\;\ldots
$$

\proclaim{Definition}
Let $\tilde{c}_0$, $ A_1$, and $w$, be any fixed real numbers and define
$\tilde{c}_j$, $j\in {\Bbb N}$, by
$\tilde{c}_j=\tilde{c}_0-j\,A_1$.

A system of real numbers
$B^n_s$, where $n$ is a non-negative integer and
$s$ is a sequence of non-negative integers with finite number
of nonzero terms,
is defined by induction with
respect to $n+[s]$ as follows
$$
\align
B^0_0=\ &1;\\
B^n_{s}=\
&(1-\de_{s_0,0})(n+|s|-1)(\tilde{c}_{n+|s|-2}-w)
\left[\sum_{l=0}^{n-1}
B^l_{s-\si_0}\right]+\\
&\sum_{i=1}^{\infty} (1-\de_{s_i,0})(s_{i-1}+1)\sum_{l=0}^{n-1}
B^l_{s-\si_i}.
\endalign
$$
In the formula above, we use the convention that
any sum $\sum_a^b...$ vanishes whenever $a>b$.
\endproclaim

In the sequel, we shall use the $B$'s with the numbers $A_1$ and $\tilde
c_0=c_0$ chosen as in \nmb!{4.4} and \nmb!{5.1}, respectively. Note
that then the numbers $B^n_s$ still depend implicitly on the value of
the variable $w$ which plays the role of the conformal weight.

The induction above works fine, because the smallest possible
value of $n+[s]$ is achieved only for $n=0$, $s=(0)$ and
the value of $B^0_0$ is fixed as $1$ in advance. The inductive
formula for $B^n_{s}$ clearly uses only
$B$'s with a smaller value of $n+[s]$.

Certain values of $B$'s are immediately clear from definition:
$B^n_{(0)}=0$ for all $n\not=0$ and $B^0_{s}=0$ for all
$s\not=(0)$.
More generally, we get from the definition by induction (with respect
to $n$) that
$B^n_s=0$ for all $n$, $s$ such that $n<[s]$.

\subhead\nmb.{7.8} Basic properties of $B^n_s$\endsubhead
Before treating more complicated examples, we shall introduce one
more piece of notation.
For a positive integer $n$,
the symbol $\{n\}$ will denote the number
$$\{n\}:=n(\tilde{c}_{n-1}-w).$$
Later on, we shall consider values of these factors $\{n\}$ at
special values of conformal weight $w=\tilde{c}_{k-1}$, $k\in\NN$.
Let us note already at this point that for this value of $w$
the resulting number depends linearly on $A_1$ (see \nmb!{7.6}).

\subhead The case where $|s|=1$\endsubhead
Using the shorthands $\{n\}$, we get immediately from the definition
that
$$\gather
B^n_{(1)}=\{n\},\quad\text{for all $n\geq 1$,}
\\
B^n_{(2)}=\{n+1\}\sum_{l=1}^{n-1}\{l\},
\quad\text{for $n\geq 2$,}
\endgather
$$
while $B^1_{(2)}=0$.

Similarly (by induction with respect to $n$), we get easily
for any $n\geq m+1$
$$
B^n_{(m+1)}=
\{n+m\}\sum_{l_{m}=m}^{n-1}\{l_m+m-1\}\!\!\sum_{l_{m-1}=m-1}^{l_m-1}\!\!
\{l_{m-1}+m-2\}\!\!\sum_{l_{m-2}=m-2}^{l_{m-1}-1}\!\ldots\sum_{l_1=1}^{l_{2}-1}\{l_1\},
$$
and $B^n_{(m+1)}=0$ for $n=0,\ldots,m.$
Clearly,
the numbers $B^n_{(m)}|_{w=\tilde{c}_{k-1}}$ are homogeneous of degree $m$
in $A_1$ for each $k\in\NN$.

\subhead The case where $|s|=2$\endsubhead
To understand the definition of $B^n_s$ better, let us also consider
the numbers $B^n_{(ij)}.$ Couples $(ij)$ of non-negative
integers can be considered as vertices of a graph in the plane.
These vertices will be connected with arrows of length $1$ going
horizontally right and antidiagonal arrows of length $\sqrt{2}$
going up and left.

Any vertex in the lattice can be reached from
$(00)$ by one or more paths (lying completely in the first quadrant).
For every path to a vertex $(ij)$, it is possible to deduce
 a  contribution to the value of $B^n_{(ij)}$
corresponding to this path
from the algorithm defining $B$'s. The actual value of
$B^n_{(ij)}$ is then the sum of such contributions over
all possible paths from $(00)$ to $(ij)$.

The situation for longer sequences $s$ is similar. The 
numbers $B^n_s$  play a principal role in the evaluation
of coefficients for standard operators, so we shall
study them in more details in Appendix B and we shall give an explicit
formula for them there.

Using the very definition of $B$'s and the simple relations
$|s-\si_0|=|s|-1$, $|s-\si_i|=|s|$, for all $i>0$,
we get immediately by induction with respect
to values of $n$ and $|s|$ the following important fact:

\proclaim{Lemma}
The numbers $B^n_s$
evaluated at $w=\tilde{c}_{k-1}$ are homogeneous of degree $|s|$
in $A_1$.
\endproclaim

\subhead\nmb.{7.9} Formulae for expansion coefficients\endsubhead
Let $k\in\NN$ be fixed.
Suppose that $j\in\NN$ and
$s=(s_0,s_1,\ldots,s_m)$ is a finite sequence of non-negative
integers
such that $j+[s]=j+\sum_{i=0}^ms_i(i+2)\leq k$.
Let $\tilde{c}_i$ be the real numbers defined in \nmb!{7.6}
and $B^n_s$ the numbers defined in \nmb!{7.7}.
Then we have the following theorem.

\proclaim{Theorem}
The coefficients $a^{k,j}_s$ in the expression  for
$D(\la,\th,k)t$
in \nmb!{7.4} are given by the formulae
$$\gather
a^{k,j}_{s}:=
\pmatrix k\\j
\endpmatrix
\left\[\prod_{i=k-j}^{k-1}({\tilde c}_i-w)\right\]
\left[\sum_{l=0}^{k-j-|s|}
B^{l}_{s}\right],\quad\text{for all $j\ge1$}            \tag"(1)"
\\
a^{k,0}_{s}:=
\sum_{l=0}^{k-|s|}
B^{l}_{s}.            \tag"(2)"
\endgather$$
\endproclaim

\demo{Proof}
The theorem will be proved by induction with respect to $k$,
using the recursive relations
from \nmb!{7.5}.

Let $k=1$. Then, according to Corollary \nmb!{7.2},
$F^1=\pi(\na t+(\tilde{c}_0-w)\ta\ox t)$.
The inequality  $j+\sum_{i=0}^ms_i(i+2)\leq 1$ is satisfied
only for $s=(0)$ and $j=0,1$.
The relations (1) and (2) read as 
$a^{1,0}_0=B_0^0+B^1_0$ and $a^{1,1}_0=(\ti{c}_0-w)B_0^0$. 
The definition of $B$'s yields $B_0^0=1,B^1_0=0$ which
proves the claim in this case.

Suppose now that the theorem holds for some fixed $k$.
Let us first prove the relation (2), i.e\. suppose first $j=0$.
By inductive assumption and the recursive relations \nmb!{7.5} for
$a$'s, we get 
$$
\multline
a^{k+1,0}_s=
\left[\sum_{l=0}^{k-|s|}B^l_{s}\right]
+(1-\delta_{s_0,0})\pmatrix
k\\1\endpmatrix({\tilde c}_{k-1}-w)
\left[\sum_{l=0}^{k-|s|}B^l_{s-\si_0}\right]
+\\
\sum_{i=1}^m (1-\delta_{s_i,0})(s_{i-1}+1)
\left[\sum_{l=0}^{k-|s|}B^l_{s-\si_i}\right]=
\sum_{l=0}^{k+1-|s|}B^l_s,
\endmultline$$
where we use
$$
\align
B^{k+1-|s|}_{s}=\ &
(1-\delta_{s_0,0})k({\tilde c}_{k-1}-w)
\left[\sum_{l=0}^{k-|s|}B^l_{s-\si_0}\right]
+\\
&\sum_{i=1}^m (1-\delta_{s_i,0})(s_{i-1}+1)
\left[\sum_{l=0}^{k-|s|}B^l_{s-\si_i}\right].
\endalign$$

For positive $j$, we get
$$
\align
a^{k+1,j}_{s}
=\ &\pmatrix   k\\j-1\endpmatrix\prod_{k-j+1}^{k-1}({\tilde c}_i-w)
\left[\sum_{l=0}^{k+1-j-|s|}B^l_{s}\right]
({\tilde c}_k-w-(k-j+1)A_1)+\\
&\;\;\;+\pmatrix   k\\j\endpmatrix\prod_{k-j}^{k-1}({\tilde c}_i-w)
\left[\sum_{l=0}^{k-j-|s|}B^l_{s}\right]
+\\
&\;\;\;+(j+1)(1-\delta_{s_0,0})\pmatrix
k\\j+1\endpmatrix\prod_{k-j-1}^{k-1}({\tilde c}_i-w)
\left[\sum_{l=0}^{k-j-|s|}B^l_{s-\si_0}\right]
+\\
&\;\;\;+\sum_{i=1}^m (1-\delta_{s_i,0})(s_{i-1}+1)
\pmatrix   k\\j\endpmatrix\prod_{k-j}^{k-1}({\tilde c}_i-w)
\left[\sum_{l=0}^{k-j-|s|}B^l_{s-\si_i}\right]
\\\allowdisplaybreak
=&\pmatrix   k+1\\j\endpmatrix\prod_{k-j+1}^{k-1}({\tilde c}_i-w)
\left[\sum_{l=0}^{k-j-|s|}B^{l}_{s}\right]
\cdot\\\allowdisplaybreak
&\;\;\;\;\;
\cdot\left[\frac{j}{k+1}({\tilde c}_k-w-(k-j+1)A_1)+\frac{k-j+1}{k+1}
({\tilde c}_{k-j}-w)
\right]+\\\allowdisplaybreak
&\;\;\;+\pmatrix   k+1\\j\endpmatrix\prod_{k-j+1}^{k-1}({\tilde c}_i-w)
\left[B^{k+1-j-|s|}_{s}\right]\cdot\\
&\;\;\;\;\;\cdot
\left[\frac{j}{k+1}({\tilde c}_k-w-(k-j+1)A_1)+\frac{k-j+1}{k+1}({\tilde c}_{k-j}-w)
\right]\\\allowdisplaybreak
=&\pmatrix   k+1\\j\endpmatrix\prod_{k-j+1}^{k}({\tilde c}_i-w)
\left[\sum_{l=0}^{k+1-j-|s|}B^{l}_{s_1,\ldots,s_m}\right],
\endalign$$
where we have used the relations
$$
\align
B^{k+1-j-|s|}_{s}=&
(1-\delta_{s_0,0})({\tilde c}_{k-j-1}-w)\sum_{l=0}^{k-j-|s|}
B^{l}_{s-\si_0}(k-j)+\\
&\sum_{i=1}^m(s_{i-1}+1)(1-\delta_{s_i,0})\sum_{l=0}^{k-j-|s|}
B^{l}_{s-\si_i}.\rlap{\qquad\qquad\qed}
\endalign$$
\enddemo

\subhead\nmb.{7.10} Formulae for the operators $D(\la,\th,k)$\endsubhead
Note that the form of the coefficients $a^{k,j}_s$ shows immediately that
all obstruction terms vanish at once for the value $w=\ti{c}_{k-1}$
of the (generalized) conformal weight. It confirms once more that
the operators $D(\la,\th,k)$ are invariant, independently of
the algebraic proof worked out in Section \nmb!{5}.
Theorem \nmb!{7.9} gives at the same time
the values of coefficients in the correction terms, i.e\. the explicit
form of the operators $D(\la,\th,k)$. It is sufficient to use
\nmb!{7.9}.(2) and to substitute there the corresponding value of $w$.

As a consequence of Lemma \nmb!{7.8} and the definition of the constants
$a^{k,0}_s$, it is clear that $a^{k,0}_s$ are homogeneous
of degree $|s|$ in $A_1$. Hence the constants $A_1$
can be absorbed into the definition of the deformation 
tensor $\Ga$ by introducing news tensors $\tilde\Ga:=A_1\Ga$
and the resulting formula is uniform and universal for all
AHS structures (for conformal structures, the constant $A_1$
is equal to $1$).

For practical calculations of curvature correction terms
of standard operators, it is better to first write
down formulas for coefficients $B^n_s$, because they
have the same form for all $k$. Having $k$ fixed,
it is then easy to evaluate $B^n_s$ at $w=\tilde{c}_{k-1}$
and to get the necessary coefficients $a^{k,0}_s$.
Note, however, that for operators of order bigger than 10,
it is better to implement the algorithm on a computer, 
since the list of correction terms is going quickly to be
unmanageable. We have postponed the exposition of the
general formulae for $B^n_s$ to Appendix B, but let
us illustrate the procedure by a few examples now.

In order to make the dependence on the order $k$ and the corresponding
fixed conformal weight $w$ explicit, we shall use the notation $B^n_s(k)$,
or $\{n\}(k)$,
for the numbers $B^n_s$, or $\{n\}$, evaluated with
$w=\tilde c_{k-1}$, respectively. Clearly
$\{n\}(k)=n(k-n)A_1$.
The numbers $B^n_{s}(k)$ are simplified
considerably, because the term $\tilde{c}_{j-1}-w$
reduces to $k-j$.
Note that after such  substitution,
'symmetric' products $\{j\}=j(k-j)A_1$
are appearing repeatedly in formulas for $B^n_s(k)$.
This leads to further simplifications of the formulae
for some $B(k)$'s, for example
$
B^n_{(n)}(2n)=[(2n-1)!!]^2
$.

\subhead\nmb.{7.11} Examples in low degrees\endsubhead
Let us recall that $B^n_s=0$ for all $n$, $s$ such that $n<[s]$ and
$B^n_{(0)}=0$ for all $n>0$.
We have already seen special cases of the previous general formulae: 
$$
B^n_{(1)}=\{n\},\;B^n_{(2)}=\{n+1\}\sum_{\ell-1}^{n-1}\{\ell\}.
$$
The Example in Appendix B provides the coefficients
$$
\gather
B^n_{(01)}= \sum_{l=1}^{n-1}\{l\}; \quad
B^n_{(001)}=\sum_{l'=2}^{n-1}\sum_{l=1}^{l'-1}\{l\}\\
B^n_{(11)}=2\sum_{l'=2}^{n-1}\{l'+1\}\sum_{l=1}^{l'-1}\{l\}+
\{n+1\} \sum_{l'=2}^{n-1}\sum_{l=1}^{l'-1}\{l\}
.\endgather$$

We denote by $\tilde\Ga$ here the corrected tensor $A_1\Ga$ and we
compute the universal formula for the operators $D(\la,\th,k)$
independently of the choice of AHS structure and the data $\la,\th$
for low values of $k$.   
The projection $\pi$ denotes as before the projection onto
the unique irreducible component $\VV_{\mu}$ in
$\ox^k(\gog_1^{\Bbb C})\ox\VV_{\la}$, the operator $D$ is
written using the conventions set up in Theorem \nmb!{7.4}, and we
write $a^k_s$ instead of $a^{k,0}_s$.
Note that by formula (2) of theorem \nmb!{7.9} we have
$a^k_{(0)}=\sum_{l=0}^kB^l_{(0)}=B^0_{(0)}=1$. 

\noindent{\bf The case $k=2$.}
Here we only need the coefficients $a^2_{(0)}=1$ and 
$$
a^2_{(1)}=B^1_{(1)}=\{1\}(2)=1.
$$
Hence
$$
D(\la,\th,2)t=\pi[\na^2t+\tilde\Ga\otimes t].
$$

\noindent{\bf The case $k=3$.}
We need the 3 coefficients $a^3_{(0)}=1$, 
$$
a^3_{(1)}=B^1_{(1)}+B^2_{(1)}=\{1\}+\{2\}\text{\ and\ } 
a^3_{(01)}=B^2_{(01)}=\{1\}.
$$
Using $\{1\}(3)=2$, $\{2\}(3)=2$, we get
$$
D(\la,\th,3)t=\pi[\na^3t+4\tilde\Ga\ox (\na t)+ 2(\na\tilde\Ga)\ox t].
$$

\noindent{\bf The case $k=4$.}
Now, we need 5 coefficients: $a^4_{(0)}=1$, and 
$$
\alignedat2
&a^4_{(1)}=B^1_{(1)}+B^2_{(1)}+B^3_{(1)}=\{1\}+\{2\}+\{3\}\quad 
&a^4_{(2)}=B^2_{(2)}=\{3\}\{1\}\\
&a^4_{(01)}=B^2_{(01)}+B^3_{(01)}=2\{1\}+\{2\}
&a^4_{(001)}=B^3_{(001)}=\{1\}.
\endalignedat
$$
Hence using again $\{n\}(k)=n(k-n)A_1$, we get
$$
D(\la,\th,4)t=\pi[\na^4t
+10\tilde\Ga\ox (\na^2 t)
+10(\na\tilde\Ga)\otimes(\na t)
+9\tilde\Ga\otimes\tilde\Ga\otimes t
+3(\na^2\tilde\Ga)\ox t].
$$

\noindent{\bf The case $k=5$.}
Here we need 7 coefficients: $a^5_{(0)}=1$, and 
$$
\align
&a^5_{(1)}=B^1_{(1)}+\ldots +B^4_{(1)}=\{1\}+\{2\}+\{3\}+\{4\}\\
&a^5_{(2)}=B^2_{(2)}+B^3_{(2)}=\{3\}\{1\}+\{4\}(\{1\}+\{2\})\\
&a^5_{(01)}=B^2_{(01)}+B^3_{(01)}+B^4_{(01)}=3\{1\}+2\{2\}+\{3\}\\
&a^5_{(001)}=B^3_{(001)}+B^4_{(001)}=\{1\}+(2\{1\}+\{2\})\\
&a^5_{(0001)}=B^4_{(0001)}=\{1\}\\
&a^5_{(11)}=B^3_{(11)}=2\{3\}\{1\}+\{4\}\{1\}
\endalign
$$
Hence we get
$$
\align
D(\la,\th,5)t=\pi[&\na^5t
+20\tilde\Ga\ox(\na^3t)
+30(\na\tilde\Ga)\ox(\na^2 t)
+64\tilde\Ga\ox\tilde\Ga\ox (\na t)
+\\
&18(\na^2\tilde\Ga)\ox(\na t)
+4(\na^3\tilde\Ga)\otimes t
+64\tilde\Ga\ox (\na\tilde\Ga)\otimes t].
\endalign
$$
As a further illustration we include the final formula in order seven. Here
we use the concatenation of the symbols instead of the tensor products and
we omit the projection $\pi$
$$\gather
\na^7 t + 
56 \tilde\Ga \na^5t +
140(\nabla\tilde\Ga) \na^4t +
168(\na^2\tilde\Ga) \na^3t +
784 (\tilde\Ga) ^2\na^3t +
2352 \tilde\Ga (\na\tilde\Ga) \na^2t +
\\
112(\na^3\tilde\Ga )\na^2t +
2304 (\tilde\Ga )^3\na t +
1180(\na\tilde\Ga) ^2\na t +
1408 \tilde\Ga (\na^2\tilde\Ga) \na^t +
40(\na^4\tilde\Ga) \na t +
\\
708(\na\tilde\Ga)(\na^2\tilde\Ga) t+
312 \tilde\Ga (\na^3\tilde\Ga ) t +
3456 (\tilde\Ga) ^2(\na\tilde\Ga)  t +
6(\na^5\tilde\Ga)  t 
\endgather$$
%
%

\head Appendix A.\endhead

For explicit description of all weights in the representation
$\gog_1$ in individual cases,
 we shall use  results gathered in \cite{FH}. The facts which
 are not proved below can be found there.

\subhead {A.1} Conformal case, even dimension\endsubhead
Here $\gog^{\CC}=\jdn$, $\semis=\dn$.
Let $L_1,\ldots,L_n$ be the standard  basis for
the dual of the Cartan subalgebra. The fundamental weights
$\pi_i,\;i=1,\ldots,n$ are
given by relations
$$
\pi_i=L_1+\ldots+L_i;\;i=1,\ldots,n-2;
\;\pi_n+\pi_{n-1}=L_1+\ldots+L_{n-1};\;
\pi_n-\pi_{n-1}=L_n.
$$
The dimension of $\gog_1$ is $2n$ and
the list of all weights of $\gog_1$
(all with multiplicity one) is given by
$\{\pm L_i;\;i=1,\ldots,n\}$.
In terms of fundamental weights, we get
$$
\align
&L_1=\pi_1;\;L_i=\pi_i-\pi_{i-1},\,i=2,\ldots,n-2;\\
&\;L_{n-1}=\pi_n+\pi_{n-1}-\pi_{n-2};
\;L_n=\pi_n-\pi_{n-1}.
\endalign
$$
Hence all  coefficients in the decompositions are in absolute
values at most one.
All weights of $\gog_1$ belong in this case to the same orbit
of the Weyl group.

\subhead {A.2} Conformal case, odd dimension\endsubhead

Here $\gog^{\CC}=\jbn$, $\semis=\bn$.
Let $L_1,\ldots,L_n$ be the standard  basis for
the dual of the Cartan subalgebra. The fundamental weights
$\pi_i,\;i=1,\ldots,n$ are
given by relations
$$
\pi_i=L_1+\ldots+L_i;\;i=1,\ldots,n-1;
\;\pi_n=(1/2)[L_1+\ldots+L_{n-1}].
$$
The dimension of $\gog_1$ is $2n+1$ and
the list of all weights of $\gog_1$
(all with multiplicity one)
 is given by
$\{0;\;\pm L_i;\;i=1,\ldots,n\}$.
In terms of fundamental weights, we get
$$
L_1=\pi_1;\;L_i=\pi_i-\pi_{i-1},\,i=2,\ldots,n-1;
\;L_n= 2\pi_n-\pi_{n-1}.
$$

So it not true in this case that all weights of $\gog_1$
have coefficients (with respect to  fundamental weights)
in absolute value less or equal to $1.$
There are two orbits of the Weyl group in the set of
all weights of $\gog_1.$ All nonzero weights form the first
orbit and the zero weight the second one.

\subhead {A.3} Grassmannian case\endsubhead
Here $\gog^{\CC}=A_{p+q+1}$, $\semis=A_p\times A_q$.
This is the only case, where $\semis$ is not a simple Lie algebra.
Irreducible representations
$V_{\la,\la'}$ of $\semis$ are just tensor
products $V_{\la}\ox V_{\la'}$ of two
irreducible representations $V_{\la}$,
resp\. $V_{\la'}$
of  $A_p$, resp\. $A_q$.
To decompose the product $V_{\la,\la'}\ox\gog_1$ means to
decompose individual products $V_{\la}\ox V$ and $V_{\la'}\ox V'$,
where $V$,  resp\. $V'$ are  defining representations of both
parts of $\semis$ and then to multiply both decompositions.

So it is sufficient to study just the case $A_n.$
Let us consider the algebra $A_n=\jan.$
Let $L_1,\ldots,L_{n+1}$ be the canonical  basis for $\CC^{n+1}.$
The dual of the Cartan subalgebra can be identified
with the quotient $\{(L_i)\in\CC^{n+1}\}/\{\sum_{i=1}^{n+1}L_i=0\}.$
The fundamental weights
$\pi_i,\;i=1,\ldots,n$ are
given by relations
$$
\pi_i=L_1+\ldots+L_i;\;i=1,\ldots,n.
$$
The dimension of the defining representation $V$ of
$\jan$ is $n+1$ and
the list of all weights of $\gog_1$
(all with multiplicity $1$)
 is given by
$\{\pm L_i;\;i=1,\ldots,n+1\}$.
In terms of fundamental weights, we get
$$
L_1=\pi_1;\;L_i=\pi_i-\pi_{i-1},\,i=2,\ldots,n;
\;L_{n+1}=-\pi_n.
$$

Hence all  coefficients in the decompositions are in absolute
values at most one.
All weights of $\gog_1$ belong in this case to the same orbit
of the Weyl group.

\subhead {A.4} Symplectic case\endsubhead
Here $\gog^{\CC}=\cn,\semis=\ban,$ hence the algebra $\semis$ is again of type
$A_k.$
Let $L_1,\ldots,L_n$ be the canonical  basis for
the defining representation $V=\CC^n$.
The dual of the Cartan subalgebra is again identified
with the quotient $\{(L_i)\in\CC^{n}\}/\{\sum_{i=1}^{n}L_i=0\}$.
The fundamental weights
$\pi_i,\;i=1,\ldots,n-1$ are
given by relations
$$
\pi_i=L_1+\ldots+L_i;\;i=1,\ldots,n-1.
$$
In this case, the representation $\gog_1$ of $\semis$ is equivalent
to $\odot^2(V)$ and its highest weight is equal to $2\pi_1$.
The dimension of $\gog_1$ is equal to $(n+1)n/2$ and
the list of all weights of $\gog_1$
(all with multiplicity $1$)
 is given by
$$
\{e_{ij}=L_i+L_j;\;i,j=1,\ldots,n;\;i\leq j\}.
$$
Using conventions $\pi_0=\pi_n=0,$ we can express $e_{ij}$ using
$\pi_j$ by
$$
e_{ij}=(\pi_i-\pi_{i-1})+(\pi_j-\pi_{j+1}),i\leq j.
$$

Hence $e_{ii}=2\pi_i-2\pi_{i-1}$
and the corresponding coefficients are $\pm 2$.
There are two orbits of the Weyl group --- $\{e_{ii}\}$
and $\{e_{ij}|i<j\}$.

\subhead {A.5} Spinorial case\endsubhead
Here $\gog^{\CC}=\dn,\semis=\ban$ and the algebra $\semis$ is again of type
$A_k.$
In this case, the representation $\gog_1$ of $\semis$ is equivalent
to $\La^2(V)$ and its highest weight is equal to the second
fundamental weight $\pi_2$.
The dimension of $\gog_1$ is equal to $n(n-1)/2$ and
the list of all weights of $\gog_1$
(all with multiplicity $1$)
 is given by
$\{e_{ij}=L_i+L_j;\;i,j=1,\ldots,n;\;i<j\}$.
Using the same conventions $\pi_0=\pi_n=0,$ we can express $e_{ij}$ using
$\pi_j$ by
$$
e_{ij}=(\pi_i-\pi_{i-1})+(\pi_j-\pi_{j+1});\;i\leq j.
$$
Hence all  coefficients in the decompositions are in absolute
values at most one.
All weights of $\gog_1$ belong in this case to the same orbit
of the Weyl group.

\subhead {A.6} $E_6$ case\endsubhead
Here $\gog^{\CC}=E_6$, $\semis=D_5$ and $\gog_1$ is one of the
basic (half)-spinor representations.
Its dimension is $16$. All weights  form one orbit of the
Weyl group and all their coefficients with respect to the
fundamental weights are in absolute value at most one.
The structure of the orbit as well as all these coefficients
can be found in \cite{Kr}.

\subhead {A.7} $E_7$ case\endsubhead
Here $\gog^{\CC}=E_7$ and $\semis=E_6$. All weights of $\gog_1$ form
one orbit of the Weyl group and all their coefficients
are in absolute value at most one (for details, see \cite{Kr}).

\head Appendix B.\endhead

To understand the definition of $B^n_s$ better, we discussed the
case of numbers $B^n_{(ij)}$ already in \nmb!{7.8}.
Couples $(ij)$ of non-negative
integers were considered as vertices of a graph in plane and
these vertices were connected with arrows of length $1$ going
horizontally right and antidiagonal arrows of length $\sqrt{2}$
going up and left.

Any vertex in the lattice can be reached from
$(00)$ by one or more paths.
For every path to a vertex $(ij)$, it is possible to deduce
its contribution to the value of $B^n_{(ij)}$
from the algorithm defining $B$'s. The actual value of
$B^n_{(ij)}$ is then the sum of such contributions over
all possible paths from $(0)$ to $(ij)$.
The situation for longer sequences $s$ is similar.
It would be possible to define a similar graph for
all sequences $s$, but it is not possible to draw it in more
general cases. We shall do the same in the language of sequences,
which also makes possible to prove an explicit formula for the values
of $B^n_s$, resp\. $B^n_s(k)$.

Let us first introduce a few additional notations.
Let ${\Cal A}$ denote the set of all finite
sequences (of a variable length)
 $J=(j_1,j_2,\ldots,j_{\alpha})$,
where $j_1=0$ and $j_2,\ldots,j_{\al}$ are non-negative integers and
put $|J|:=\al$.
For a positive integer $a$ and $J\in{\Cal A}$,
let us define the sequences $s^J$, $ s^J_a$ by
$$
s^J:=\sum_{a'=1}^{|J|}\si_{j_{a'}};\qquad
s^J_a:=\sum_{a'=1}^{a}\si_{j_{a'}};\;a=1,\ldots,|J|-1;\qquad s^J_0:=(0)
$$
where $\si_i$ are the sequences from \nmb!{7.7}.
The subset ${\Cal A}_0$ of ${\Cal A}$ is defined by
$$
{\Cal A}_0:=\{J\in{\Cal A}\ |\ (s^J_a)_i\geq 0;\;
a=1,\ldots,|J|,\,i=0,1,\ldots\;\}.
$$

We have the following simple properties
$$\gather
[\si_i]=1\text{ for all }i\text{ and }
[\si_i]+[\si_j]=[\si_i+\si_j]\text{ for all }i,j
\\
[s^J]=|J|.
\endgather
$$

In order to generalize formulas for $B^n_{(m)}$ deduced in Section \nmb!{7},
let us introduce for every sequence $s$ of non-negative integers
the set
$$
{\Cal A}^0_s:=\{J\in {\Cal A}^0\ |\ s^J=s\}.
$$
This set is a generalization of the set of all different paths from
$(0)$ to $s$ discussed above in the case of sequences of length two.

We also need a generalization of the numbers $\{n\}$ from \nmb!{7.8}.
Let us define the numbers
$\{s,l,a\}$, where $s$ is a finite sequence of integers
and $l$, $a$ are positive integers
$$
\{s,l,a\}:=\cases \{l+|s|\} \quad &\text{if $a=0$}\\
s_{a-1}\quad &\text{if $a\ne0$.}
\endcases
$$
Using all this notation we obtain the following explicit formula for
the numbers $B^n_s$:

\proclaim{Theorem} The numbers $B^n_s$ are given by the formula
$$
\sum_{J\in{\Cal A}^0_s}
\{s^J_{\al-1},n,j_{\al}\}
\sum_{l_{\al-1}=\al-1}^{n-1}
\{s^J_{\al-2},l_{\al-1},j_{\al-1}\}
\sum_{l_{\al-2}=\al-2}^{l_{\al-1}-1}
\ldots
\sum_{l_2=2}^{l_3-1}
\{s^J_{1},l_{2},j_2\}
\sum_{l_1=1}^{l_{2}-1}
\{l_1\}
$$
where $\al=[s]=|J|$.
\endproclaim

\demo{Proof}
We can use induction with respect to $\al$.
The case $\al=1$ means that $s=(1)$.
This case was discussed in \nmb!{7.8}:
$B^n_{(1)}=\{n\}$.
But $s=\si_0,$ there is just one element $J=(0)$ in
${\Cal A}^0_s$ and the theorem holds.

Suppose now that the formula is valid for all $s$ with $[s]\leq k-1$
and consider a sequence $s$ with $[s]=k.$
The set ${\Cal A}^0_s$ of sequences $J$  can be split into a disjoint union
of subsets by an additional condition $j_{[s]}=i,i=0,1,\ldots,$
(all but a finite number of them being empty).
Now, let us have a look at
the algorithm defining $B$'s.
Using the induction assumption for terms
$\sum_{l=0}^{n-1}B^l_{s-\si_i},\,i=0,1,\ldots$
and noticing that
$n+|s|-1=n+|s-\si_0|;\;
s_{i-1}+1=(s-\si_i)_{i-1},
$
we get the correct value for $B^n_s.$
\hfill\qed
\enddemo

\subhead Examples \endsubhead Let us use the formula in a few cases.
If $s=(01)$, then the set ${\Cal A}^0_s$ is a one point set. It consists
of $J=(0,1)$, $s=\si_0+\si_1$. Hence
$$
B^n_{(01)}=\{(1),n,1\}\sum_{l=1}^{n-1}\{l\}=\sum_{l=1}^{n-1}\{l\}.
$$
Similarly, for $s=(001)$, we have
${\Cal A}^0_s=\{(0,1,2)\}$, $s=\si_0+\si_1+\si_2$.
Hence
$$
B^n_{(001)}=
\{(01),n,2\}\sum_{l'=2}^{n-1}\{(1),l',1\}\sum_{l=1}^{l'-1}\{l\}=
\sum_{l'=2}^{n-1}\sum_{l=1}^{l'-1}\{l\}.
$$

If $s=(11)$, there are two elements in the set ${\Cal A}^0_s,$
namely $J=(0,0,1),\,s=\si_0+\si_0+\si_1$
and $J=(0,1,0),\,s=\si_0+\si_1+\si_0.$ So
$$
\align
B^n_{(11)}&=
\{(2),n,1\}\sum_{l'=2}^{n-1}\{(1),l',0\}\sum_{l=1}^{l'-1}\{l\}+
\{(01),n,0\}\sum_{l'=2}^{n-1}\{(1),l',1\}\sum_{l=1}^{l'-1}\{l\}=\\
&=2\sum_{l'=2}^{n-1}\{l'+1\}\sum_{l=1}^{l'-1}\{l\}+
\{n+1\} \sum_{l'=2}^{n-1}\sum_{l=1}^{l'-1}\{l\}.
\endalign
$$

A similar computation leads to the last constant $B^4_{(0001)}=\{1\}$ which
we have used in \nmb!{7.11}.

\enddocument